\documentclass[a4paper,11pt]{article}
\usepackage{amssymb}
\usepackage{amsfonts}
\usepackage{amsmath,tabu, color}
\usepackage{amsthm}
\usepackage{cite}

\newcommand{\s}{s}
\newcommand{\eps}{{\varepsilon}}

\newcommand{\R}{\mathbb{R}}

\newcommand{\N}{\mathbb{N}}

\newcommand{\cuad}{{\sqcap\kern-.68em\sqcup}}

\newcommand{\norm}[1]{\|#1\|}

\numberwithin{equation}{section}
\newtheorem{theorem}{Theorem}[section]
\newtheorem{definition}[theorem]{Definition}
\newtheorem{proposition}[theorem]{Proposition}

\newtheorem{lemma}[theorem]{Lemma}
\newtheorem{corollary}[theorem]{Corollary}
\newtheorem{remark}[theorem]{Remark}
\newcommand{\bremark}{\begin{remark} \em}
\newcommand{\eremark}{\end{remark} }

\newcommand{\cA}{{\mathcal A}}

\newcommand{\cE}{{\mathcal E}}
\newcommand{\cF}{{\mathcal F}}

\newcommand{\cH}{{\mathbb H}}

\newcommand{\cL}{{\mathcal L}}
\newcommand{\cM}{{\mathcal M}}

\begin{document}
\begin{center}{\bf  \Large On the bounds of the sum of eigenvalues for a Dirichlet problem involving mixed fractional Laplacians}\medskip

\bigskip\medskip

 {\small    Huyuan Chen\footnote{chenhuyuan@yeah.net}
\smallskip

   School of Mathematics and Statistics, Jiangxi Normal University, Nanchang,\\
Jiangxi 330022, PR China \\[4mm]

  Mousomi Bhakta\footnote{ mousomi@iiserpune.ac.in} \smallskip

  Department of Mathematics, Indian Institute of Science Education and Research (IISER-Pune),\\
Pune 411008, India \\[4mm]

 Hichem Hajaiej\footnote{  hichem.hajaiej@gmail.com}
\smallskip

California State University, Los Angeles, 5151, USA\\[6mm]

}

\medskip

\begin{abstract}
Our purpose in this paper is to study of the eigenvalues $\{\lambda_i(\mu)\}_i$ of the  Dirichlet problem
$$ (-\Delta)^{s_1} u=\lambda \big( (-\Delta)^{s_2}  u+\mu u\big)\ \  {\rm in}\ \,  \Omega,\quad\quad  u=0\ \ {\rm in}\  \ \R^N\setminus \Omega,$$
where $0<s_2<s_1<1$, $N>2s_1$ and $(-\Delta)^{s}$ is the fractional Laplacian operator defined in the principle value sense.

 We first show the existence of a sequence of eigenvalues, which approaches infinity. Secondly we provide a Berezin--Li--Yau  type  lower bound for the   sum of the eigenvalues of the above problem. Furthermore, using a self-contained and novel method, we establish an  upper bound for the sum of eigenvalues of the problem under study.

 \end{abstract}
  \end{center}

 \vspace{1mm}


\vspace{1mm}

\setcounter{equation}{0}
\section{Introduction and main results}

Let $0<\s_2<\s_1<1$, $\Omega$ be  a smooth bounded domain in $\R^N$ with the integer $N\ge 1$, $N>2s_1$. The main goal of   this paper is to study
 the lower bounds of the eigenvalues for the Dirichlet  problem
\begin{equation}\label{eq 1.1}
\left\{ \arraycolsep=1pt
\begin{array}{lll}
 (-\Delta)^{\s_1} u=\lambda  \Big( (-\Delta)^{\s_2} u+\mu u\Big)\quad \  &{\rm in}\quad   \Omega,\\[2mm]
 \phantom{(-\Delta)^\s   }
  u=0\quad \ &{\rm{in}}\  \quad \R^N\setminus \Omega,
\end{array}
\right.
\end{equation}
where  $(-\Delta)^s$ is the fractional laplacian   defined in the following sense  (principle value)
 \begin{equation}\label{fl 1}
 (-\Delta)^\s  u(x)=c_{N,\s} \lim_{\epsilon\to0^+} \int_{\R^N\setminus B_\epsilon(x) }\frac{u(x)-
u(y)}{|x-y|^{N+2\s}}  dy
\end{equation}
with $s\in(0,1)$,
$c_{N,\s}=2^{2\s}\pi^{-\frac N2}\s\frac{\Gamma(\frac{N+2\s}2)}{\Gamma(1-\s)}$ and
$\Gamma$  being the  Gamma function, see e.g. \cite{RS1}.
 Recall that, for $\s\in(0,1)$, the fractional Laplacian of a function $u \in C^\infty_c(\R^N)$ can  also  be defined by:
\begin{equation*}
  \label{eq:Fourier-representation}
(-\Delta)^\s u(\xi) =\mathcal{F}^{-1}\bigg( |\xi|^{2\s}\mathcal{F}(u) (\xi)\bigg)\qquad \text{for all $\xi \in \R^N$},
\end{equation*}
where $F(f)$ denotes the Fourier transform of $f$.

\medskip

 Problem  (\ref{eq 1.1}) involves two fractional Laplacians with two different powers.
The terminology  ``mixed operators" refers to the differential or pseudo-differential order of the operator,
and to the type of the operator, which can combine classical and fractional features.
 When $\lambda<0$, (\ref{eq 1.1}) involves a sum of two fractional Laplacians of different orders. Indeed,  such sum of operators  arises naturally from superposition of two stochastic processes with a different random walk and a L\'evy flight, this is the case when a particle can follow either of these two processes according to a certain probability, the associated limit diffusion equation is described by a sum of  fractional Laplacians, see for e.g \cite{B-1} .    While if $\lambda>0$, (\ref{eq 1.1}) models a difference of  two fractional Laplacians. The sum and the difference of two fractional Laplacians appear in many circumstances.  To mention few, problems of blood circulation in the heart, responsible for causing heart problems and in many circumstances coronary bypass surgeries, can be modelled by two to five mixed fractional Laplacians, for e.g. see \cite{E-3, M-4,M-5,M-6}. It also appears in many circumstances because of the anomalous blood circulation, but it is often not the same anomaly in all the arteries, and the blood can follow either of the five arteries. Other applications of mixed fractional operators with different orders include plasma physics and population dynamics, ways to reduce pandemics and so on. In view of these important applications, we strongly believe that equation (\ref{eq 1.1}) and some of its variants described above will get an increasing interest in the near future.

The  most primitive model of (\ref{eq 1.1})  for the eigenvalues     is
\begin{equation}\label{eq 1.1-s=1}
\left\{ \arraycolsep=1pt
\begin{array}{lll}
 -\Delta u=\lambda u \quad \  &{\rm in}\quad   \Omega,\\[1.5mm]
 \phantom{ -\Delta   }
  u=0\quad \ &{\rm{on}}  \quad\! \partial\Omega,
\end{array}
\right.
\end{equation}
which has attracted the attention of mathematicians  since 1912.  Indeed, {  Weyl in \cite{W}  showed  that  the $k$-th eigenvalue, $\lambda_{1,k}(\Omega)$ ($1$ stands for $s=1$) of \eqref{eq 1.1-s=1}, }has the asymptotic behavior
$$\lambda_{1,k} (\Omega)\sim C_N(k|\Omega|)^{\frac2N}\quad{\rm as}\ \  k\to+\infty,$$
 where
 $$C_N=(2\pi)^2|B_1|^{-\frac{2}{N}}$$
 and $|B_1|$ is the volume of unit ball in $\R^N$. Later, P\'olya \cite{P} (in 1960) provided  a lower bound for the k-th eigenvalue,
\begin{equation}\label{p-conj}
\lambda_{1,k}(\Omega)\geq C(k/|\Omega|)^{\frac2N}
\end{equation}
for any { {\it plane covering domain}} $D$ in $\R^2$ with $C=C_N$, (his proof also works in dimension $N\geq 3$). $D$ is called a { {\it plane covering domain} in $\R^2$ if an infinity of domains congruent to $D$ cover the whole plane without gaps and without overlapping except a set of measure zero. In \cite{P}, P\'olya  also made a conjecture that (\ref{p-conj}) holds for any bounded domain in $\R^N$ with $C=C_N$.}  To answer this conjecture,  Lieb \cite{L} proved (\ref{p-conj}) with a positive constant $C$ for general bounded domain and Li-Yau \cite{LY} improved
 the constant $C$ to $\frac{N}{N+2}C_N$.  Until now, this constant for lower bound is the best  and (\ref{p-conj}) with $C=\frac{N}{N+2}C_N$  is   called  Berezin-Lieb-Yau inequality.   More related estimates on lower bounds for the eigenvalues under various setting can be found in \cite{CgW,F1,L,M}.  On the other hand, the upper bounds of Dirichlet eigenvalues were derived by Kr\"oger in \cite{K} by calculating the Rayleigh quotient by using a sequence of functions approaching the characterized function of $\Omega$.    We also refer to Yang's upper bounds of the Dirichlet's eigenvalues in  \cite{CgY,CQL} in the following way:
 $$\lambda_{1,k}(\Omega)\leq c(N,k)  \lambda_{1,1}(\Omega)k^\frac{2}{N} \quad\text{  for  some $c(N,k)>0$}.$$

  The unstopped interest in finding bounds for eigenvalues of the Dirichlet problem is in part due to the following fact: The Hilbert-Polya conjecture is to associate the zero of the Riemann Zeta function with the
eigenvalue of a Hermitian operator. This quest initiated the mathematical interest for estimating
the sum of Dirichlet eigenvalues of the Laplacian while in physics the question is related to count
the number of bound states of a one body Schr\"odinger operator and to study their asymptotic distribution. The latter constitutes in itself a branch in nonlinear analysis.
During the last  decade, there has been a renewed and increasing
interest in the study of linear and nonlinear integral operators.
The prototype is the fractional Laplacian. This has been motivated by numerous applications in different fields motivated by connections to real world life applications and by important advances in the theory of linear and nonlinear
partial differential equations, see basic properties \cite{Hajaiej1, Hajaiej2, Hajaiej3, musina-nazarov}, regularities   \cite{CS0,RS}, Liouville property   \cite{C}, general nonlocal operator \cite{CT},
fractional Pohozaev identity \cite{RS1}, singularities \cite{BCP, CFQ,CQ}, uniqueness \cite{Frank}, fractional variational setting \cite{BM-2, EGE,F2,JXL,SV} and the references therein.

 When $\mu=0$, $s_2=0$ and $s=s_1\in(0,1)$, (\ref{eq 1.1}) reduces to
the fractional  Laplacian  problem
\begin{equation}\label{eq 1.1-s}
\left\{ \arraycolsep=1pt
\begin{array}{lll}
 (-\Delta)^{\s} u=\lambda u \quad \  &{\rm in}\quad   \Omega,\\[2mm]
 \phantom{(-\Delta)^\s   }
  u=0\quad \ &{\rm{in}}\  \quad \R^N\setminus \Omega,
\end{array}
\right.
\end{equation}
  for which the asymptotic behavior of  eigenvalues $\lambda_{s,i}$
 has been studied,    for Klein-Gordon operator i.e. $s=\frac12$ in \cite[Proposition 3.1]{HY} or for general order $s\in(0,1)$ in \cite[Theorem 1]{G},
 \begin{equation}\label{asym frac}
\lim_{k\to+\infty} k^{ -\frac{2s}{N}} \lambda_{s,k}  = a(N,s)  |\Omega|^{-\frac{2s}{ N}}
\end{equation}
where
 \begin{equation}\label{asym frac-c}
 a(N,s)= (2\pi)^{2s}|B_1|^{{-\frac{2s}{N}}}.
 \end{equation}
Moreover, a refinement of Berezin--Li--Yau-type lower bound for the sum of eigenvalues was built by Yolcu and Yolcu in \cite[Theorem 1.4]{YY} as follows
 \begin{equation}\label{low frac}
 \sum^{k}_{j=1}\lambda_{s,j}  \geq \frac{N}{N+2s} a(N,s) |\Omega|^{-\frac{2s}{ N}}  k^{ 1+\frac{2s}{N}} +ck^{1-\frac{2-2s}{N}}
\end{equation}
for some $c>0$ depending on $|\Omega|$. In a recent work \cite{YH},  Hajaiej and Wang  provided the asymptotic behavior of  the sum of the eigenvalues of  (\ref{eq 1.1-s})
  \begin{equation}\label{low frac-1}
  \lim_{k\to+\infty}k^{-1-\frac{2s}{N}} \sum^k_{j=1}\lambda_{s,j}(\Omega)=\frac{N}{N+2s} a(N,s)   |\Omega|^{-\frac{2s}{N}}.
  \end{equation}
For more estimates on eigenvalues of the fractional Dirichlet problem,  we refer the readers to \cite{CZ,G,HY,YY,YY1}, a review \cite{Frank} and the references therein.

To analyze  the fractional Dirichlet eigenvalues of (\ref{eq 1.1}),  we denote $\cH^s_0(\Omega)$   the space of all measurable functions $u:\R^N\to \R$  with $u \equiv 0$ in $\R^N \setminus \Omega$ and
$$
{ \iint _{\R^{2N}}} \frac{|u(x)-u(y)|^2}{|x-y|^{N+2s}} dx dy <+\infty.
$$
{  It is well known that} $\cH^s_0(\Omega)$ is a Hilbert space  equipped with the inner product
$$  \mathcal{E}_s(u,w)=\frac{c_{N,s}}2 { \iint _{\R^{2N}}} \frac{(u(x)-u(y)) {(w(x)-w(y))}}{|x-y|^{N+2s}}
dx dy $$
and the induced norm $$\norm{u}_s:=\sqrt{\mathcal{E}_s(u,u)}.$$  

We say a function $u \in \cH^{s_1}_0(\Omega)$  be an eigenfunction of (\ref{eq 1.1}) corresponding to the eigenvalue
$\lambda$ if
 $$ \cE_{s_1}(u,w)  = \lambda \Big( \cE_{s_2}(u, w)  + \mu \int_{\Omega} uw\,dx \Big) \quad \text{for all $w\in \cH^{s_1}_0(\Omega)$.} $$

 For $\mu\geq -\lambda_{s_2,1}$, {  we denote by  $\cH^{s_2}_{\mu,0}(\Omega)$ the space of all measurable functions $u:\R^N\to \R$  with $u \equiv 0$ in $\R^N \setminus \Omega$ and
 $$\|u\|_{s_2,\mu}:=\bigg(\cE_{s_2}(u, u)  + \mu \int_{\Omega} u^2\,dx\bigg)^\frac{1}{2}<\infty.
 $$
The corresponding inner product in $\cH^{s_2}_{\mu,0}(\Omega)$ is given by}
 $${ \langle u, w\rangle_{s_2,\mu}:=}\,  \cE_{s_2}(u, w)  + \mu \int_{\Omega} uw\,dx,\quad \forall\,\, u,\, w\in  \cH^{s_2}_0(\Omega).$$
Let $\lambda_{s_2,1}$ be the first eigenvalue of
\eqref{eq 1.1-s} corresponding to $s=s_2$.
 {  We note that $\|.\|_{s_2, \mu}$ is equivalent to $\|.\|_{s_2}$ for $\mu>-\lambda_{s_2,1}$.}

\smallskip

 Our first aim is to show the existence of a sequence of discrete  eigenvalues of (\ref{eq 1.1}) as follows.

\begin{theorem}\label{teo 1}
 Let $\mu> -\lambda_{s_2,1}$, where $\lambda_{s_2,1}>0$ be the first eigenvalue of (\ref{eq 1.1-s}) with $s=s_2$.
Then problem (\ref{eq 1.1}) admits a sequence of real eigenvalues
$$
0<\lambda_{1} (\mu)\leq \lambda_{ 2} (\mu)\le \cdots\le \lambda_{j} (\mu)\le \lambda_{ j+1} (\mu)\le \cdots
$$
and the corresponding eigenfunction $\phi_i$, $i \in \N$. Moreover, we have the following properties:
\begin{enumerate}
\item[(i)] $\lambda_{j}(\mu)=\min \{\cE_{s_1}(u,u) :\,  u\in \cH_{0,j}(\Omega), \,  {  \norm{u}_{\cH^{s_2}_{\mu, 0}(\Omega)}=1}\}$, where
$$
{  \cH_{0,1}(\Omega)=\cH^{s_1}_0(\Omega), \quad  \cH_{0,j}(\Omega):=\{u\in\cH^{s_1}_0(\Omega)\::\: \langle u, \,  \phi_m \rangle_{s_2, \mu} =0 \, \, \mbox{for}\, \,  m=1,\dots, j-1\}\,\, \text{for}\,\, j>1;}
$$
\item[(ii)] $\{\phi_j\::\: j \in \N\}$ is an orthonormal basis of $\cH^{s_2}_{\mu,0}(\Omega)$;
\item[(iii)] $\lim \limits_{j\to \infty} \lambda_{j} (\mu)=+\infty$;
\item[(iv)]  For $\mu\in(-\lambda_{s_2,1},+\infty)$,  the map $\mu\mapsto \lambda_1(\mu)$ is decreasing and $\displaystyle\limsup_{\mu\to { -\lambda_{s_2,1}^+}} \lambda_{1} (\mu)<+\infty$.
\end{enumerate}
 \end{theorem}

 We remark that
\begin{enumerate}
\item[(a)] from the appendix in \cite{SV},  problem (\ref{eq 1.1-s})  has the property that the first eigenvalue is simple and
 the corresponding eigenfunction $\phi_{s,1}$ is positive; these  properties  are derived by the following:  $\cE_{s}(|\phi_{s,1}|)< \cE_{s}(\phi_{s,1})$, if $\phi_{s,1}$ is sign-changing.
 However, this argument fails for problem (\ref{eq 1.1}) due to presence of  multiple fractional Laplacians, and  it is  very interesting but challenging to determine the one-fold of $\lambda_{1} (\mu)$  and the positivity of the eigenfunction  $\phi_1$ corresponding to the first eigenvalue $\lambda_1$  for problem (\ref{eq 1.1});
 \item[(b)]  thanks to the monotonicity and boundedness of $\lambda_{1} (\mu)$, assertion $(iv)$  indicates
that it is possible to obtain the existence of $\{\lambda_{j} (\mu)\}_{j\in\N}$ for $\mu\leq -\lambda_{s_2,1}$;

 \item[(c)] it is known that eigenfunctions of  (\ref{eq 1.1-s}) are $C^\infty(\Omega)$. To see this, one uses bootstraps method to prove solutions of \eqref{eq 1.1-s} are in $L^\infty(\Omega)$ and then uses regularity results of   \cite{RS}. While the regularity for (\ref{eq 1.1}) seems to be much more complicated,  because  bootstraps iteration has to work between different order fractional Laplacians.
\end{enumerate}

 We now provide a lower bound for the sum of eigenvalues of (\ref{eq 1.1}).

\begin{theorem}\label{teo 1.1}
Let $\mu\geq 0$ and  $\{\lambda_{j}(\mu)\}_{j\in\N}$  be the increasing sequence of eigenvalues of  problem (\ref{eq 1.1}) and $\omega_{_{N-1}}$ denote the surface area of the unit sphere in $\R^N$.  Then   for $k\in\N$
\begin{align} \label{low b}
 \sum^k_{j=1}\lambda_j(\mu)   \geq b_1  |\Omega|^{- \frac{2(s_1-s_2)}{N}}  k^{1+\frac{2(s_1-s_2)} {2\s_2+N}} - \mu   b_2  |\Omega|^{- \frac{2s_1}{ N}}  k^{1+ \frac{2\s_1-4s_2 } {2\s_2+N}},
\end{align}
where
$$b_1 =  \frac{(2\s_2+N)^{\frac{2\s_1 +N} {2\s_2+N}}}{2\s_1 +N } \Big(2^{-(N+1+2s_2)}\pi^{-\frac{3N}2}\omega_{_{N-1}}^{2-\frac{2s_2}N}     \frac{\Gamma(\frac{N-2s_2}{2})}{\Gamma(s_2+1)}    N^{\frac{2s_2}N}\Big)^{-\frac{2(\s_1 -s_2)} {2\s_2+N}}$$
and
$$b_2  =b_1  \frac{2s_1+N}{N(2s_2+N)^{\frac{2s_2}{2s_2+N}}}  \Big(2^{-(N+1+2s_2)}\pi^{-\frac{3N}2}\omega_{_{N-1}}^{2-\frac{2s_2}N}     \frac{\Gamma(\frac{N-2s_2}{2})}{\Gamma(s_2+1)}    N^{\frac{2s_2}N}\Big)^{\frac{2s_2}{N+2s_2}}. $$
\end{theorem}

\medskip

We remark that when $s_2=0=\mu$, the constant $b_1(s_1, 0)$ in Theorem~\ref{teo 1.1} reduces to  $\frac{N}{N+2s_1}a(N, s_1)$, where $a(N,s_1)$ is defined as in \eqref{asym frac-c}. Moreover, \eqref{low b} coincides with Berezin-Li-Yau estimate for \eqref{eq 1.1-s=1} (with $s_1=1$).

The proof of Theorem~\ref{teo 1.1} is inspired by the Berezin--Li--Yau method
(see \cite{LY}) in which the authors treated mainly the function $\Psi_k(x,y)=\sum^k_{j=1}\psi_j(x)\psi_j(y), $
where $\psi_j$'s are the eigenfunctions corresponding to the eigenvalues $\lambda_j$ of \eqref{eq 1.1-s=1}. Denote $F=\int_{\R^N} |(\cF_x\Psi_k)(z,y)|^2 dy$,
 the key estimate is the following
\begin{equation}\label{LY-key}\int_{\R^N} Fdz\leq \Big(\frac{N+2}{N} \int_{\R^N}|z|^2 Fdz\Big)^{\frac{N}{N+2}} \Big(|B_1|\|F\|_{L^\infty}\Big)^{\frac{2}{N+2}},
\end{equation}
where $|B_1|$ denotes the volume of unit ball in $\R^N$. For our problem (\ref{eq 1.1}),  the situation is much more complicated. More precisely, to apply Berezin-Li-Yau method,
 we need to consider the function $f=\int_{\R^N} |(\cF_x\Phi_k)(z,y)|^2 dy$, where
 $$\Phi_k(x,y)=\sum^k_{j=1}\phi_j(x)\tilde \phi_j(y), \quad \tilde \phi_j:= \Big((-\Delta)^{s_2}+\mu\Big)^{\frac12}\phi_j$$  and $\phi_j$'s are eigenfunctions corresponding to the eigenvalues $\lambda_j$ of \eqref{eq 1.1}. The most important estimate is
 $\int_{\R^N} (|z|^{2s_2}+\mu) f(z)dz$, which is controlled by $\int_{\R^N}|z|^{2s_1} fdz$ and $\|f\|_{L^\infty}$.
These difficulties arise from  the non homogeneous property of lower oder term $(|z|^{2s_2}+\mu) f$ and the estimate $\|f\|_{L^\infty}$.   For  $\|f\|_{L^\infty}$,  the original tool is the Bessel's inequality, which requires  orthonormal property for $\{\phi_j\}$ in $L^2(\Omega)$, while $\{\phi_j\}_{j}$ is not an orthonormal sequence in $L^2(\Omega)$  (it is orthonormal in $H^{s_2}_{\mu,0}(\Omega)$). To overcome this difficulty, we transform $\phi_j$ to $\tilde \phi_j$ which is orthonormal in $L^2(\R^N)$ but not supported in $\Omega$. Further using certain delicate estimate, we have overcome the difficulty (see Section~3.2).

\medskip

The following is a direct corollary from Theorem \ref{teo 1.1} using the monotonicity of map $j\to \lambda_{j}(\mu)$.

 \begin{corollary}\label{cr 1.1}
Let $\mu\geq 0$ and  $\{\lambda_{j}(\mu)\}_{j\in\N}$  be the increasing sequence of eigenvalues of  problem (\ref{eq 1.1}).   Then   for $k\in\N$
\begin{align} \label{low b-k}
\lambda_k (\mu)   \geq b_1  |\Omega|^{- \frac{2(s_1-s_2)}{ N}}  k^{\frac{2(s_1-s_2)} {2\s_2+N}} - \mu   b_2  |\Omega|^{- \frac{2s_1}{ N}}  k^{ \frac{2\s_1-4s_2 } {2\s_2+N}},
\end{align}
where $b_1$ and $b_2$ are same as in Theorem~\ref{teo 1.1}.
 \end{corollary}

For the upper bounds, due to the numerous challenges mentioned earlier, we only address the case   that $\mu=0$ and  a restriction on the upper range of  $s_1$. More specifically, we have the following estimates on upper bounds.

  \begin{theorem}\label{teo upp}
Let $0<s_2<s_1<  \frac{1+s_2}2  $, $\mu=0$ and  $\{\lambda_{i}(\mu)\}_{i\in\N}$  be the increasing sequence of eigenvalues of  problem (\ref{eq 1.1}) and $\Omega$ be a bounded $C^2$ domain.
Then there exists $c_0=c_0(N,s_1, s_2, \Omega)>0$ and $\delta_3\in(0, 1+\frac32\frac{ \s_1-s_2} { N})$ such that
    for $k\in\N$
\begin{align} \label{upp b}
 \sum^k_{j=1}\lambda_j (0)   \leq b_3 |\Omega|^{-\frac{2(s_1-s_2)}{N}}  k^{1+\frac{2(\s_1-s_2)} { N}} + c_0  k^{\delta_3},
\end{align}
where
$$b_3 = (2\pi)^{2(s_1-s_2)} \omega_{_{N-1}}^{-\frac{2(s_1-s_2)}{N}}\frac{N^{1+\frac{2(s_1-s_2)}{N}}}{N+2(s_1-s_2)}. $$
 \end{theorem}

It is worth noting that
 \begin{enumerate}
\item[(d)]  the constant $b_3$ in the upper bound (\ref{upp b})  coincides   with  \eqref{low frac-1}  with $s=s_1-s_2$;
\item[(e)]  the upper bound and the lower bound for our problem (\ref{eq 1.1}) obtained in Theorem \ref{teo 1.1} and Theorem  \ref{teo upp}  are not enough to determine $ \sum^k_{j=1}\lambda_j (\mu)$, even for $\mu=0$.

 \item[(f)]  from our proofs,  it isn't too difficult to see that Theorem~\ref{teo 1} and Theorem~\ref{teo 1.1} can be extended to the case $s_1=1$.  However,  it fails for the upper bound in Theorem \ref{teo upp}
 because of the restriction of $s_1<  \frac{1+s_2}2<1$, which is not essential since it appears  due to the technique difficulty.
 \end{enumerate}
 
 \medskip
 
The rest of the paper is organized as follows. In Section 2, we  study  qualitative properties of the eigenvalues and the corresponding eigenfunctions  of problem (\ref{eq 1.1}), namely Theorem~\ref{teo 1}. Section 3 is devoted to show the lower bound of the sum of eigenvalues, namely Theorem~\ref{teo 1.1}. Finally, in section 4, we discuss the upper bounds aspects, namely Theorem~\ref{teo upp}.

\medskip

{\bf Notations}: Throughout this paper, $\omega_{N-1}$ denotes the surface area of unit sphere in $\R^N$, $B_r(x) \subset \R^N$ is an open ball of radius $r$ centered at $x \in \R^N$, and we set $B_r:= B_r(0)$ for $r>0$. For any set $A$ of $\R^N$, $|A|$ denotes the Lebesgue measure of $A$ and $\mathcal{F}(f)$ denotes the Fourier transform of a function $f$. By $C^\infty_c(\R^N)$, we denote $C^\infty$ functions in $\R^N$ with compact support.

\section{Existence  }

We set $Q:=\R^{2N}\setminus (\Omega^c \times \Omega^c)$, where $\Omega^c=\R^N \setminus \Omega$ and for $1\leq p<\infty$ define
$$
{ W^{s,p}_0(\Omega)}:=\Big\{u:\mathbb{R}^N\to\mathbb{R}\mbox{ measurable }\Big| u=0\ {\rm a.e.\ in\ } \Omega^c \mbox{ and }
\int_{Q} \frac{|u(x)-u(y)|^p}{|x-y|^{N+sp}}dxdy<\infty\Big\}.
$$
Note that from the fractional Poincar\'e inequality, see \cite{EGE},  the space ${ W^{s,p}_0(\Omega)}$ is endowed with the norm defined as
$$\|u\|_{W^{s,p}_0(\Omega)}:= \left(\int_{Q} \frac{|u(x)-u(y)|^p}{|x-y|^{N+sp}}dxdy\right)^{1/p}.$$

\begin{theorem}\label{teo 1.2}
Let $0<\s_2<\s_1<1$, $p>1$ and $\Omega$ be  a bounded domain in $\R^N$, then ${  W^{s_1,p}_0(\Omega) \hookrightarrow W^{s_2,p}_0(\Omega)}$ is continuous and compact.

\end{theorem}

Before starting to prove the above theorem, we need to introduce the Besov Space over $\Omega$. We follow the notations of \cite[Section 2.3.1]{Tr}.

\begin{definition}
Let $\mathcal{S}$ denote the Schwartz class functions on $\R^N$ and 
 $\cA $ be collection of all sequence $\eta=\{\eta_i\}_{i=0}^\infty\in \mathcal{S}(\R^N)$ such that
\begin{center}
supp$(\eta_0)\subset \{x: |x|\leq 2\}$,    supp$(\eta_j)\subset \{x: 2^{j-1}\leq |x|\leq 2^{j+1}\}$\ \ for $j=1,2, \cdots$,
\end{center}
$$\sum_{i=0}^{\infty}\eta_j(x)=1,$$
and for every multi-index $\alpha$ there exists a positive number $c_\alpha$ such that
$$2^{j\alpha}|D^\alpha\eta_j(x)|\leq c_\alpha \quad\forall\, j=0,1,2,\dots \quad\mbox{and}\quad\forall\, x\in\R^N.$$
\end{definition}

\begin{definition}
Let $s\in(-\infty, \infty)$ and $p,q\in(0,\infty]$ and $\{\eta_i\}_{i=0}^\infty\in \cA $. Then
$$B^s_{p,q}(\R^N):=\Big\{f\in \mathcal{S'}(\R^N): \|f\|_{B^s_{p,q}(\R^N)}^\eta:=\|2^{sj}\cF^{-1}\big(\eta_j\cF(f)\big)\|_{l_q\big(L^p(\R^N)\big)}<\infty\Big\}.$$
\end{definition}
It can be shown that the quasi-norm {  $\|f\|_{B^s_{p,q}(\R^N)}^\eta$} does not depend on the choice of {  $\eta\in\cA$} (see \cite[Proposition 1 in 2.3.2]{Tr}).

\begin{definition} Let $\mathcal{D'}(\Omega)$ denote the set of all distributions over $\Omega$. For $1\leq p,\, q \leq \infty$, and $0<s<1$,  we set
$$B^s_{p,q}(\Omega)= \Big\{u \in \mathcal{D'}(\Omega) : \,\exists\,\, g \in B^s_{p,q}(\R^N) \quad\mbox{with}\quad g|_{\Omega}=u \Big\}$$
and $$\|u\|_{B^s_{p,q}(\Omega)}=\inf_{g \in B^s_{p,q}(\R^N),\, g|_{\Omega}=u} \|g\|_{B^s_{p,q}(\R^N)}.$$	
Here $B^s_{p,q}(\Omega)$ is called the Besov Space over $\Omega$.
\end{definition}	

\begin{lemma} \label{l:1} \cite[Theorem 3.1.1(i)]{Tr}
Let $p_0,q_0, p_1, q_1\in (0, \infty)$, $s_0, s_1\in (-\infty, \infty)$ and $\Omega$ be a smooth bounded domain in $\R^N$. Then the following embedding
$$B^{s_0}_{p_0,q_0}(\Omega)\hookrightarrow B^{s_1}_{p_1,q_1}(\Omega)$$
is continuous if $s_0-\frac{N}{p_0}>s_1-\frac{N}{p_1}$.
\end{lemma}

\begin{lemma} \label{l:2} \cite[pg.233]{Tr}
Let $0<s_2<s_1$ and $p,q\in (0, \infty)$ and $\Omega$ be a smooth bounded domain in $\R^N$. Then the embedding $B^{s_1}_{p,q}(\Omega)\hookrightarrow B^{s_2}_{p,q}(\Omega)$ is compact.
\end{lemma}

 \noindent {\bf Proof of Theorem \ref{teo 1.2}. } From \cite[Lemma 2.2]{BM}, we know that the embedding ${  W^{s_1, p}_0(\Omega)\subset W^{s_2, p}_0(\Omega)}$ is continuous. Moreover, from \cite[pg. 209]{Tr}), it follows that
 $$\|u\|_{W^{s,p}(\Omega)} :=\|u\|_{L^p(\Omega)}+\displaystyle\bigg(\int_{\Omega\times\Omega}\frac{|u(x)-u(y)|^p}{|x-y|^{N+sp}}dxdy \bigg)^{\frac{1}{p}}$$
is an equivalent norm for $\|u\|_{B^s_{p,p}(\Omega)}$ for $p\in(1,\infty), \,s\in(0,1)$.
Therefore, by Lemma~\ref{l:2} we have that
\begin{equation}\label{em:3}
W^{s_1, p}(\Omega)\hookrightarrow W^{s_2, p}(\Omega) \quad\text{is compact.}
\end{equation}
 Now let $\{u_n\}$ be a bounded sequence in ${ W^{s_1, p}_0(\Omega)}$, to prove the theorem we need to extract a convergent subsequence in
${  W^{s_2, p}_0(\Omega)}$. By Rellich compactness, up to a subsequence $u_n\to u$ in $L^p(\Omega)$ for some $u\in { W^{s_1, p}_0(\Omega)}$. For that subsequence we define $v_n:= u_n-u$. Therefore,
 we need to show $v_n\to 0$ in ${  W^{s_2, p}_0(\Omega)}$. As $v_n=0$ in $\Omega^c$,
 \begin{eqnarray}\label{eq:4}
 \|v_n\|_{{  W^{s_2, p}_0(\Omega)}} ^p &=&\int_{\Omega\times\Omega}\frac{|v_n(x)-v_n(y)|^p}{|x-y|^{N+s_2p}}dxdy+2\int_{x\in\Omega}\int_{y\in\Omega^c}\frac{|v_n(x)-v_n(y)|^p}{|x-y|^{N+s_2p}}dxdy\nonumber\\
 &\leq &\|v_n\|_{W^{s_2, p}(\Omega)}^p+2\int_{x\in\Omega}\int_{y\in\Omega^c}\frac{|v_n(x)-v_n(y)|^p}{|x-y|^{N+s_2p}}dxdy.
 \end{eqnarray}

{  By P\'oincare inequality in \cite[Theorem 6.7]{EGE} there exists a positive constant $c>0$ such that}
  $$\|v\|_{L^p(\Omega)}^p\leq c \int_{\Omega\times\Omega}\frac{|v(x)-v(y)|^p}{|x-y|^{N+sp}}dxdy\quad {\rm for\ any}
   \ v\in {  W^{s, p}_0(\Omega)}.$$
 Thus, from the definition of ${  W^{s_1, p}_0(\Omega)}$ and $W^{s_1, p}(\Omega)$, there exists $c>0$ such that
$$\|v_n\|_{W^{s_1, p}(\Omega)}\leq c\|v_n\|_{{ W^{s_1, p}_0(\Omega)}}.$$
 Therefore, $\{v_n\}$ is a bounded sequence in $W^{s_1, p}(\Omega)$. Consequently,
 by \eqref{em:3}, up to a subsequence $v_n\to 0$ in $W^{s_2, p}(\Omega)$. Therefore,
 to complete the proof, we only have to show that
 $$\int_{x\in\Omega}\int_{y\in\Omega^c}\frac{|v_n(x)-v_n(y)|^p}{|x-y|^{N+s_2p}}dxdy\to 0.$$ To that end, let $\epsilon>0$ be arbitrary and we define $M:=\sup_n \|v_n\|_{{ W^{s_1, p}_0(\Omega)}}$. Then
 \begin{eqnarray*}
 \int_{x\in\Omega}\int_{y\in\Omega^c}\frac{|v_n(x)-v_n(y)|^p}{|x-y|^{N+s_2p}}dxdy&=& \int_{x\in\Omega}\int_{y\in\Omega^c\cap |x-y|<\eps}\frac{|v_n(x)-v_n(y)|^p}{|x-y|^{N+s_2p}}dxdy \\
 &&\quad+ \int_{x\in\Omega}\int_{y\in\Omega^c\cap |x-y|\geq\epsilon}\frac{|v_n(x)-v_n(y)|^p}{|x-y|^{N+s_2p}}dxdy\\
 &=&\int_{x\in\Omega}\int_{y\in\Omega^c\cap |x-y|<\epsilon}\frac{|v_n(x)-v_n(y)|^p}{|x-y|^{N+s_1p}}|x-y|^{(s_1-s_2)p}dxdy\\
 &&\quad+ \int_{x\in\Omega}\bigg(\int_{y\in\Omega^c\cap |x-y|\geq\epsilon}\frac{dy}{|x-y|^{N+s_2p}}\bigg)|v_n(x)|^pdx\\
 &\leq& \epsilon^{(s_1-s_2)p}\iint_{\R^N\times\R^N}\frac{|v_n(x)-v_n(y)|^p}{|x-y|^{N+s_1p}}dxdy \\
 &&\quad+C(\epsilon)\|v_n\|_{L^p(\Omega)}^p\\
  &\leq&\epsilon^{(s_1-s_2)p}M+o(1),
 \end{eqnarray*}
 where $o(1)\to 0$ as $n\to\infty$. Therefore, the above integral can be made arbitrary small. Hence from \eqref{eq:4}, we conclude that $ \|v_n\|_{{ W^{s_2, p}_0(\Omega)}}\to 0$. This completes the proof. \hfill$\Box$

\medskip

In the particular case that $p=2$, for $s\in(0,1)$ we set
$$\cH^{s}_0(\Omega)={  W^{s, 2}_0(\Omega)}$$
which is a Hilbert space with the inner product
$\cE_{s}(u,v)$ for $u,v\in \cH^{s}_0(\Omega)$.

\medskip

\noindent{\bf Proof of Theorem \ref{teo 1}. }
The functional $$
\cH^{s_1}_0(\Omega) \to \R, \qquad u \mapsto \Phi(u):= \cE_{s_1}(u,u)
$$
is weakly lower semicontinuous. Moreover, let
$$\cM_1:= \{u\in \cH^{s_1}_0(\Omega),\, \norm{u}_{s_2} + \mu \norm{u}_{L^2(\Omega)}=1\},$$
where $\norm{\cdot}_{s_2} + \mu \norm{\cdot}_{L^2(\Omega)}$ is equivalent to
$\norm{\cdot }_{s_2}$ in the space $\cH^{s_2}_0(\Omega)$ for $\mu>-\lambda_{s_2,1}$.
Then we have that
$$
 \Phi(u)  < +\infty \qquad \text{for $u \in \cM_1$.}
$$
Put $\lambda_1(\mu) := \inf_{\cM_1} \Phi$.  { By Theorem \ref{teo 1.2}, the embedding
$W^{s_1, p}_0(\Omega)\hookrightarrow W^{s_2, p}_0(\Omega)$
is compact and therefore}, it follows that $\lambda_1(\mu) $ is attained by a function $\phi_1 \in \cM_1$. Consequently, there exists a Lagrange multiplier $\lambda \in \R$ such that
$$
\cE_{s_1}(\phi_1, w)=\frac{1}{2}\Phi'(\phi_1)w= \lambda \Big( \cE_{s_2}(\phi_1, w)  + \mu \int_{\Omega} \phi_1 w\,dx \Big)  \qquad \text{for all $w\in \cH^{s_1}_0(\Omega)$.}
$$
Choosing $w = \phi_1$ yields to $\lambda= \lambda_1(\mu) $. Hence $\phi_1$ is an eigenfunction of (\ref{eq 1.1}) corresponding to the eigenvalue $\lambda_1 (\mu)$.   Moreover $\lambda_1(\mu)>0$. Next we proceed inductively and assume that $\phi_2,\dots,\phi_k \in \cH^{s_1}_0(\Omega)$ and $\lambda_2 (\mu)  \le \dots \le \lambda_k(\mu)  $ are already given for some $k \in \N$ with the properties that for $j=2,\dots,k$, the function $\phi_j$ is a minimizer
of $\Phi$ within the set
 \begin{align*}
\cM_j:= \{u\in \cH^{s_1}_0(\Omega)\::\: \norm{u}_{s_2} + \norm{u}_{L^2(\Omega)}=1,\: \text{$\cE_{s_2}(u, \phi_m)  + \mu \int_{\Omega} u\phi_m\,dx  =0$ for $m=1,\dots j-1$} \},
 \end{align*}
$\lambda_j(\mu) = \inf_{\cM_j} \Phi = \Phi(\phi_j)$, and
\begin{equation}
  \label{eq:inductive-eigenvalue}
  \cE_{s_1}(\phi_j, \varphi)=\lambda_j(\mu) \Big(\cE_{s_2}(\phi_j, \varphi)  + \mu \int_{\Omega} \phi_j \varphi \,dx\Big)    \qquad \text{for all $\varphi\in \mathbb{H}^{s_1}_0(\Omega)$.}
\end{equation}
{ Again by the compact embedding in Theorem \ref{teo 1.2}}, the value $\lambda_{k+1} $ is attained by a function { $\phi_{k+1} \in \cM_{k+1}$.} Consequently, there exists a Lagrange multiplier $\lambda \in \R$ such that
\begin{equation}
  \label{eq:inductive-eigenvalue-k+1}
\cE_{s_1}(\phi_{k+1}, \varphi)=\lambda \Big(\cE_{s_2}(\phi_{k+1}, \varphi)  + \mu \int_{\Omega} \phi_{k+1} \varphi \,dx\Big)   \qquad \text{for all $\varphi \in \cM_{k+1}(\Omega)$.}
\end{equation}
Choosing $\varphi = \phi_{k+1}$, we have that  $\lambda= \lambda_{k+1}(\mu) $. Moreover, for $j=1,\dots,k$, we have, by (\ref{eq:inductive-eigenvalue}) and the definition of $\cM_{k+1}(\Omega)$, that
\begin{align*}
\cE_{s_1}(\phi_{k+1}, \phi_j)=\cE_{s_1}(\phi_j, \phi_{k+1})&= \lambda_j(\mu)  \Big(\cE_{s_2}(\phi_j, \phi_{k+1})  + \mu \int_{\Omega} \phi_j \phi_{k+1}\,dx \Big)
\\& = 0
\\&= \lambda_{k+1} (\mu) \Big(\cE_{s_2}(\phi_{k+1}, \phi_j )  + \mu \int_{\Omega} \phi_j \phi_{k+1}\,dx \Big).
\end{align*}
Hence (\ref{eq:inductive-eigenvalue-k+1}) holds with $\lambda= \lambda_{k+1}(\mu)  $ for all $\varphi\in \cH^{s_1}_0(\Omega)$.
Inductively, we have now constructed a normalized sequence $(\phi_k)_k$ in $\cH^{s_1}_0(\Omega)$ and a nondecreasing sequence $\{\lambda_k\}_k$ in $\R$ such that property (i) holds and such that $\phi_k$ is an eigenfunction of (\ref{eq 1.1}) corresponding to $\lambda = \lambda_k(\mu)  $ for every $k \in \N$. Moreover, by construction, the sequence $\{\phi_k\}_k$ forms an orthonormal system in $\cH^{s_2}_{\mu,0}(\Omega)$.\smallskip

Next we show property (iii), i.e., $\lim \limits_{k\to+\infty} \lambda_k(\mu)  =+\infty.$ Supposing by contradiction that $c:= \lim \limits_{k \to \infty}\lambda_k(\mu)  <+\infty$, we deduce that
$\cE_{s_1}(\phi_k,\phi_k) \le c$ for every $k \in \N$. Hence the sequence $(\phi_k)$ is bounded in $\cH^{s_1}_0(\Omega)$, and therefore by Rellich compactness theorem, $(\phi_k)$
contains a convergent subsequence $(\phi_{k_j})_j$ in $L^2(\Omega)$. This however is impossible since the functions $\{\phi_{k_j}\}_{j \in \N}$ are orthonormal in $H^{s_2}_{\mu,0}(\Omega)$. Hence (iii) is proved.\smallskip

Next, to prove that $\{\phi_k\::\: k \in \N\}$ is an orthonormal basis of $H^{s_2}_{\mu,0}(\Omega)$, we first suppose by contradiction that there exists $v \in \mathbb{H}^{s_1}_0(\Omega)$ with $\|v\|_{s_2}+{ \mu}\|v\|_{L^2(\Omega)} = 1$
 and $\cE_{s_2}(v, \phi_k )  + \mu \int_{\Omega} v \phi_k\,dx   =0$ for any $k \in \N$.
Since $\lim \limits_{k\to\infty} \lambda_k (\mu)=+\infty$, there exists an integer $k_0>0$  such that
$$\Phi(v)<\lambda_{k_0} (\mu)=\inf_{\cM_{k_0}}\Phi(u),$$
 which by definition of $\cM_{k_0}$ implies that
$\cE_{s_2}(v, \phi_k )  + \mu \int_{\Omega} v \phi_k\,dx  \not = 0$ for some $k \in \{1,\dots,k_0-1\}$. This is a a contradiction to $\cE_{s_2}(v, \phi_k )  + \mu \int_{\Omega} v \phi_k\,dx   =0$ for any $k \in \N$. Thus, we conclude that $\cH^{s_1}_0(\Omega)$ is contained in the $H^{s_2}_{\mu,0}$-closure of the span of $\{\phi_k\::\: k \in \N\}$. Since $\cH^{s_1}_0(\Omega)$ is dense in $H^{s_2}_{\mu,0}(\Omega)$, we conclude that the span of $\{\phi_k\::\: k \in \N\}$ is dense in $H^{s_2}_{\mu,0}(\Omega)$, and hence $\{\phi_k\::\: k \in \N\}$ is an orthonormal basis of $H^{s_2}_{\mu,0}(\Omega)$. This proves (ii). 

\smallskip

Finally,  we show $(iv)$. It follows by the definition of the first eigenvalue $\lambda_1(\mu)$ 
that for $\mu\in(-\lambda_{s_2,1},+\infty)$, the map $\mu\mapsto \lambda_1(\mu)$ is decreasing.

Let $(\lambda_{s_2,1}, \varphi_{s_2,1})$ be the first eigenvalue and the corresponding eigenfunction  of (\ref{eq 1.1-s}) with $s=s_2$.  We may assume that $\varphi_{s_2,1}>0$ in $\Omega$,  $\int_\Omega\varphi_{s_2,1}^2dx=1$.   Let $\eta_0$ be a smooth function with compact support in $B_{r_0}(x_0)$ such that $B_{2r_0}(x_0)\subset \Omega$
with $\int_\Omega\eta_0^2dx=1$. Note that $\eta_0\not \equiv \varphi_{s_2,1}$ in $\Omega$.
Then $\cE_{s_2}(\eta_0,\eta_0)>\lambda_{s_2,1}$, otherwise $\eta_0\equiv\varphi_{s_2,1}$
since the first eigenfunction of   (\ref{eq 1.1-s}) with $s=s_2$   is simple.

Thus, we have that
 \begin{align*} \lambda_1(\mu)&\leq \frac{\cE_{s_1}(\eta_0,\eta_0)}{\cE_{s_2}(\eta_0,\eta_0)+\mu\int_\Omega \eta_0^2 dx}
 \\&\  \leq   \frac{\cE_{s_1}(\eta_0,\eta_0)}{\cE_{s_2}(\eta_0,\eta_0)-\lambda_{s_2,1} }
\  <+\infty \quad{\rm as}\ \ \mu\to -\lambda_{s_2,1}^+.
 \end{align*}
This completes the proof.
\hfill$\Box$

 \section{Lower bounds}\label{ }
\subsection{Important estimate}

To prove the Li-Yau's type lower bound for \eqref{eq 1.1}, we need the following results.

\begin{proposition}\label{pr 2.1}
Let $\mu\geq0$, $0<\s_2<\s_1<1$ and  $f$ be a real-valued function defined on $\R^N$ with $0\leq f\leq M_1$,
$$\int_{\R^N} f(z)|z|^{2s_1} dz\leq M_2,$$
 then
 \begin{align*}
 \int_{\R^N} (|z|^{2s_2}+\mu) f(z)dz  & \leq  \frac{\omega_{_{N-1}} M_1}{2\s_2+N}   \Big(\frac{(2\s_1+N) M_2}{  M_1 \omega_{_{N-1}}}\Big)^{\frac{2\s_2+N}{2\s_1 +N}} \Big[ 1+ \mu\frac{2\s_2+N}N\Big(\frac{(2\s_1+N) M_2}{  M_1 \omega_{_{N-1}}}\Big)^{-\frac{2\s_2}{2\s_1+N}}\Big].
 \end{align*}

\end{proposition}

\noindent{\bf Proof of Proposition \ref{pr 2.1}. } Let
 $$h(z)=\left\{ \arraycolsep=1pt
\begin{array}{lll}
(|z|^{2s_2}+\mu)  M_1\quad \quad  &{\rm for}\quad   |z|< R,\\[1.5mm]
 \phantom{    }
  0\quad \ &{\rm{for}}\quad   |z|\geq R,
\end{array}
\right.$$
where $R>0$ such that
$$ \int_{B_R}  M_1|z|^{2s_1} dz=M_2.$$
Direct computation shows that
$$
R^{2\s_1+N} = \frac{(2\s_1+N) M_2}{  M_1 \omega_{_{N-1}}}.
$$

  Since
 \begin{align*}
  |z|^{2s_2}  \Big(|z|^{2(s_1-s_2)}-R^{2(s_1-s_2)}\Big) \Big( f(z)-(|z|^{2s_2}+\mu)^{-1}h(z)\Big)\geq 0
 \end{align*}
 and
 \begin{align*}
 \mu  \Big(|z|^{2s_1 }-R^{2s_1}\Big) \Big( f(z)-(|z|^{2s_2}+\mu)^{-1}h(z)\Big)\geq 0,
 \end{align*}
it follows that
 \begin{align}\label{18-11-1}
& \int_{\R^N}|z|^{2s_2}\Big(  f(z) -(|z|^{2s_2}+\mu)^{-1} h(z)\Big) dz \nonumber
\\\leq&  \frac1{R^{2(s_1-s_2)}}\int_{\R^N}  |z|^{2s_2}   |z|^{2(s_1-s_2)}\Big(  f(z) -(|z|^{2s_2}+\mu)^{-1} h(z)\Big) dz\nonumber
\\\leq&  \frac1{R^{2(s_1-s_2)}} \int_{\R^N}  |z|^{2s_1} (f(z) -M_1\chi_{_{B_R}})dz\leq 0,
\end{align}
and
\begin{align}\label{18-11-2}
& \mu \int_{\R^N} \Big(  f(z) -(|z|^{2s_2}+\mu)^{-1} h(z)\Big) dz \nonumber
\\\leq& \frac{\mu}{R^{2s_1 }} \int_{\R^N}   |z|^{2s_1}  \Big(  f(z) -(|z|^{2s_2}+\mu)^{-1} h(z)\Big) dz\nonumber
\\\leq&  \frac{\mu}{R^{2s_1 }} \int_{\R^N}  |z|^{2s_1} (f(z) -M_1\chi_{_{B_R}})dz\leq 0.
\end{align}

Combining \eqref{18-11-1} and \eqref{18-11-2}, we have
 \begin{align*}
 \int_{\R^N} (|z|^{2s_2}+\mu) f(z)dz&\leq \int_{\R^N} h(z)dz\\
&= \frac{\omega_{_{N-1}} M_1}{2\s_2+N} R^{2\s_2+N} \Big(1+\mu\frac{2\s_2+N}N  R^{-2s_2}\Big)\\
&=  \frac{\omega_{_{N-1}} M_1}{2\s_2+N}   \Big(\frac{(2\s_1 +N) M_2}{  M_1 \omega_{_{N-1}}}\Big)^{\frac{2\s_2+N}{2\s_1 +N}}
 \Big[ 1+ \mu\frac{2\s_2+N}N \Big(\frac{(2\s_1 +N) M_2}{  M_1 \omega_{_{N-1}}}\Big)^{-\frac{2\s_2}{2\s_1 +N}}\Big].
 \end{align*}
We complete the proof.\hfill$\Box$

 \medskip

We also need the following Lemma.

\begin{lemma}\label{lm 2.1}
Assume that $1>\tau_1>\tau_2>0$  and   $d_1>0.$

 Let     $r_1$ be the solution of
 $$r^{\tau_1}\Big(1+ r^{-\tau_2}\Big)=d_1,$$
then
 $$  d_1^{\frac1{\tau_1}}\Big(1-\frac{ 1}{\tau_1}  d_1^{-\frac{\tau_2}{\tau_1}}\Big) <r_1<d_1^{\frac1{\tau_1}}.  $$

 \end{lemma}
\noindent{\bf Proof. } Let
$$f(r)=r^{\tau_1}\Big(1+  r^{-\tau_2}\Big),$$
then $f(0)=0$, $\lim_{r\to+\infty}f(r)=+\infty$ and  $f'(r)=\tau_1 r^{\tau_1-1}+ (\tau_1-\tau_2) r^{\tau_1-\tau_2-1}>0$, $f''(r)=\tau_1(\tau_1-1) r^{\tau_1-2}+ (\tau_1-\tau_2)(\tau_1-\tau_2-1) r^{\tau_1-\tau_2-2}<0$,  As a consequence,  $f$ is strictly increasing, concave in $(0,+\infty)$
and for any $d_1>0$, there exists a unique solution $r_1$ such that $f(r_1)=d_1$. Moreover, since $f(r)>0$ for $r>0$, we conclude $r_1>0$.

Let
$$ R_1= d_1^{\frac1{\tau_1}}\Big(1-\frac{1}{\tau_1}  d_1^{-\frac{\tau_2}{\tau_1}}\Big)_+\quad {\rm and} \quad R_2= d_1^{\frac1{\tau_1}},$$
 where $a_+=\max\{0, a\}$.
Note that
\begin{align*}
\frac{ f(R_2)}{d_1}&=  1 +  d_1^{-\frac{\tau_2}{\tau_1}} > 1
 \end{align*}
 and if
 $\frac{1}{\tau_1}  d_1^{-\frac{\tau_2}{\tau_1}} <1$,  then
 \begin{align*}
\frac{ f(R_1)}{d_1}&=  \Big(1-\frac{1}{\tau_1}  d_1^{-\frac{\tau_2}{\tau_1}}\Big)^{\tau_1}  +  d_1^{-\frac{\tau_2}{\tau_1}}\Big(1-\frac{1}{\tau_1}  d_1^{-\frac{\tau_2}{\tau_1}}\Big)^{\tau_1-\tau_2}
\\&\leq 1-   d_1^{-\frac{\tau_2}{\tau_1}}  +  d_1^{-\frac{\tau_2}{\tau_1}} -{ \frac{\tau_1-\tau_2}{\tau_1}d_1^{-\frac{2\tau_2}{\tau_1}}}\\
&<1,
 \end{align*}
  where we used the fact that for $\tau\in(0,1)$,
 $$(1-t)^\tau\leq 1-\tau t\quad {\rm for\ any}\ t\in(0,1).$$
On the other hand, if  $\frac{1}{\tau_1}  d_1^{-\frac{\tau_2}{\tau_1}} \geq1$,  then $R_1=0$. Therefore, $\frac{ f(R_1)}{d_1}=0<1$. Thus, in both the cases
$f(R_1)<d_1<f(R_2)$. Now using the fact that $f$ is strictly increasing, continuous and $f(r_1)=d_1$, we conclude $R_1<r_1<R_2$. This completes the proof. 
 \hfill$\Box$

 \subsection{Lower bound}

The Bessel inequality plays an important role in the Li-Yau method for the lower bound of eigenvalues:
{\it Let $\cH$ be an inner product space or a Hermitian product space together with its product function $\langle \cdot,\cdot\rangle$.  Let $e_1, e_2, \cdots$  be any (finite or infinite) orthonormal sequence. Then for any $x\in\cH$,
 $$\sum_{j\geq 1} |\langle x,e_j\rangle|^2 \leq \langle x,x\rangle=\|x\|_{_\cH}^2.$$} \medskip

\noindent {\bf Proof of Theorem \ref{teo 1.1}.}
Let
 $$\Phi_k(x,y)=\sum^k_{j=1}\phi_j(x)     \tilde \phi_j(y), $$
 where
  $$\tilde \phi_j:= \Big((-\Delta)^{s_2}+\mu\Big)^{\frac12}\phi_j=\cF^{-1}\bigg[\big(|z|^{2s_2}+\mu\big)^{\frac12}\cF\big(\phi_j\big)\bigg].$$

Using the Fourier transform, we have that
 $$\cF_x(\Phi_k)(z,y)=(2\pi)^{-\frac N2} \int_{ \R^N} \Phi_k(x,y)e^{ix\cdot z}dx $$
and
$$\cF\left(\Big((-\Delta)^{s_2}+\mu\Big)_x^{\frac12}{ \phi_j}\right)(z)= \big(|z|^{2s_2}+\mu\big)^{\frac12}{ \cF\big(\phi_j(z)\big)}\ \ \text{ for} \ \, \mu\geq0.$$
Therefore, doing a straightforward computation we have
 \begin{align}\label{Li-Yau con1}
   \int_{\R^N}\int_{\R^N}\Big|\big(|z|^{s_2}+\mu\big)^{\frac12}   \cF_x(\Phi_k)(z,y) \Big|^2 dy dz= \sum_{j=1}^k \int_{\R^N}\int_{\R^N} \tilde \phi_j^2(x)\tilde \phi_j^2(y)    dy dx=k,
 \end{align}
 by the orthonormality of  $\{\tilde \phi_j\}_{j\in\N}$ in $L^2(\R^N)$. Next, we   estimate $\big|\int_{\R^N}    \Big|  \cF_x(\Phi_k)(z,y)\Big|^2 dy\big|_{L^\infty(\R^N)}$. Here the main difficulty comes from the fact that $\{\phi_j\}_{j}$ is not an orthonormal sequence in $L^2(\Omega)$.

 \medskip

Let $G_\mu$ be the fundamental solution of $(-\Delta)^{s_2}+\mu$ in $\R^N$. Then
$$0<G_\mu(x)\leq G_0(x)=a_{N, s_2} |x|^{2s_2-N}\quad{\rm for}\ \, \mu\geq0, $$

\begin{equation}\label{a_Ns}a_{N,s_2}=2^{-2s_2} \pi^{-\frac N2} \frac{\Gamma(\frac{N-2s_2}{2})}{\Gamma(s_2)}=  2^{-2s_2} \pi^{-\frac N2} \frac{\Gamma(\frac{N-2s_2}{2})}{\Gamma(s_2+1)}s_2:=\tilde a_{N,s_2}s_2.\end{equation}
Note that
\begin{align*}
\big((-\Delta)^{s_2}+\mu\big)^{-1} g=G_\mu\ast g\quad{\rm in}\quad \R^N
 \end{align*}
 for $g\in L^2(\R^N)$.    Using Fourier transform, we can see that    $$\cF\bigg(\big((-\Delta)^{s_2}+\mu\big)^{-1} g\bigg)(z)=\big(|z|^{2s_2}+\mu\big)^{-1}\cF(g(z)).$$
  Consequently,  $\big((-\Delta)^{s_2}+\mu\big)^{-1/2} g$ can be defined as follows
 $$\big((-\Delta)^{s_2}+\mu\big)^{-1/2} g:=\bigg[\big((-\Delta)^{s_2}+\mu\big)^{-1}\bigg]^{1/2} g=
 \cF^{-1}\bigg[\big((|z|^{2s_2}+\mu)^{-1}\big)^{1/2}\cF(g)\bigg].$$

Now, in view of Parseval's relation and  Bessel's inequality, we obtain
\begin{align*}
  \int_{\R^N}    \Big|  \cF_x(\Phi_k)(z,y)\Big|^2 dy &= (2\pi)^{-  N} \sum^k_{j=1}\Big| \int_{\R^N}   \phi_j (x) \big( e^{{\rm i} x\cdot z}\chi_{\Omega}(x) \big) dx \Big|^2\int_{\R^N}|\tilde\phi_j(y)|^2 dy
  \\&= (2\pi)^{-  N}\sum^k_{j=1}\Big| \int_{\R^N}  \Big((-\Delta)^{s_2}+\mu\Big)_x^{\frac12} \phi_j (x) \Big((-\Delta)^{s_2}+\mu\Big)_x^{-\frac12}\big(e^{{\rm i} x\cdot z} \chi_{\Omega}(x) \big) dx \Big|^2
  \\&= (2\pi)^{-  N}\sum^k_{j=1}\Big| \int_{\R^N} \tilde \phi_j (x) \Big((-\Delta)^{s_2}+\mu\Big)_x^{-\frac12}\big(e^{{\rm i} x\cdot z} \chi_{\Omega}(x) \big) dx \Big|^2
\\  &\leq (2\pi)^{-  N}   \int_{\R^N}  \Big| \Big((-\Delta)^{s_2}+\mu\Big)_x^{-\frac12}\big(e^{{\rm i} x\cdot z} \chi_{\Omega}(x) \big)  \Big|^2 dx
\\  &= (2\pi)^{-  N}   \int_{\R^N}\bigg[ \Big((-\Delta)^{s_2}+\mu\Big)_x^{-1} \big(e^{{\rm i} x\cdot z} \chi_{\Omega}(x) \big) \bigg] \overline{\big(e^{{\rm i} x\cdot z} \chi_{\Omega}(x) \big)}    dx
\\  &\leq (2\pi)^{-  N}   \int_{\Omega} \Big|G_\mu\ast \big(e^{{\rm i} x\cdot z} \chi_{\Omega}(x) \big)   \Big|    dx
\\&\leq (2\pi)^{-  N} \int_{\Omega} G_0\ast \chi_{\Omega} dx
\\&  = (2\pi)^{-  N}  a_{N,s_2} \int_\Omega \int_\Omega |x-y|^{2s_2-N} dy dx
\\&\leq  (2\pi)^{-  N} a_{N,s_2}   |\Omega| \sup_{x\in\Omega}  \int_{\Omega }|x-y|^{2s_2-N} dy.
 \end{align*}
Now we choose $r>0$ such that $|\Omega|=\frac{\omega_{N-1}}{N}r^N$ and using rearrangement inequality, we have
\begin{align*} \sup_{x\in\Omega}  \int_{\Omega }|x-y|^{2s_2-N} dy&\leq \int_{B_r(x) }|x-y|^{2s_2-N} dy
\\&= \frac{\omega_{_{N-1}}}{2s_2} r^{2s_2} =  \frac{\omega_{_{N-1}}}{2s_2}\bigg(\frac{N|\Omega|}{\omega_{N-1}}\bigg)^\frac{2s_2}{N}=a_0 |\Omega|^{ \frac{2s_2}{N}},
\end{align*}
where
\begin{equation}\label{a_0}a_0=\frac{1}{2s_2}   N^{\frac{2s_2}N}    \omega_{_{N-1}}^{1-\frac{2s_2}N}.\end{equation}

 As a consequence, we obtain that
 \begin{align}\label{Li-Yau con2}
   \int_{\R^N}    \Big|  \cF_x(\Phi_k)(z,y)\Big|^2 dy  \leq (2\pi)^{-N} a_{N,s_2}a_0 |\Omega|^{1+\frac{2s_2}{N}}=\tilde a_{N,s_2} \frac{N^\frac{2s_2}N\omega_{_{N-1}}^{1-\frac{2s_2}N}}{2^{N+1}\pi^N   } |\Omega|^{1+\frac{2s_2}{N}}.
    \end{align}

Meanwhile, using the definition of $\Phi_k$, it also follows that
 \begin{align*}
   \int_{\R^N}\int_{\R^N} |z|^{2s_1} \big|\cF_x(\Phi_k)(z,y)\big|^2 dydz
  &= \int_\Omega \int_\Omega \Phi_k(x,y) (-\Delta)^{\s_1}_x  \Phi_k(x,y) dydx
 \\&=\sum^k_{j=1}\int_\Omega \phi_j(-\Delta)^{s_1}\phi_jdx\int_\Omega |\tilde \phi_j (y)|^2 dy
\\&   =\sum^k_{j=1}\lambda_j(\mu).
 \end{align*}

 Now we apply Proposition \ref{pr 2.1} to the function
 $$f(z)=\int_{\R^N} |(\cF_x\Phi_k)(z,y)|^2 dy $$
 with
 $$  M_1=(2\pi)^{-N}a_{N,s_2} a_0 |\Omega|^{ \frac{N+2s_2}{N}}\quad{\rm and}\quad M_2=\sum^k_{j=1}\lambda_j(\mu),$$
  then we conclude that
 \begin{align}\label{li-yau-estimate 1}
 k  & \leq  \frac{\omega_{_{N-1}} M_1}{2\s_2+N}   \Big(\frac{(2\s_1 +N) M_2}{  M_1 \omega_{_{N-1}}}\Big)^{\frac{2\s_2+N}{2\s_1 +N}} \cdot \left( 1+ \mu\frac{2\s_2+N}{N}\Big(\frac{(2\s_1 +N) M_2}{  M_1 \omega_{_{N-1}}}\Big)^{-\frac{2\s_2}{2\s_1 +N}}\right).
 \end{align}
 {\it Case: $\mu=0$. }  In this case,  (\ref{li-yau-estimate 1})   reduces to
 $$k    \leq  \frac{\omega_{_{N-1}} M_1}{2\s_2+N}   \Big(\frac{(2\s_1 +N) M_2}{  M_1 \omega_{_{N-1}}}\Big)^{\frac{2\s_2+N}{2\s_1 +N}},$$
 that is
$$M_2\geq  \frac{(2\s_2+N)^{\frac{2\s_1 +N} {2\s_2+N}}}{2\s_1 +N }    (M_1 \omega_{_{N-1}})^{-\frac{2(\s_1 -s_2)} {2\s_2+N}}    k ^{\frac{2\s_1 +N} {2\s_2+N}},$$
 which implies that
 \begin{align*}
 \sum^k_{j=1}\lambda_j(0)&\geq
 \frac{(2\s_2+N)^{\frac{2\s_1 +N} {2\s_2+N}}}{2\s_1 +N } \big((2\pi)^{-N} a_{N, s_2} a_0\omega_{_{N-1}} \big)^{-\frac{2(\s_1 -s_2)} {2\s_2+N}}    |\Omega|^{-\frac{N+2s_2}{N}\frac{2(s_1-s_2)}{N+2s_2}} k^{1+\frac{2(s_1-s_2)} {2\s_2+N}}
 \\&=b_1 |\Omega|^{- \frac{2(s_1-s_2)}{ N}} k^{1+\frac{2(s_1-s_2)} {2\s_1+N}}.
 \end{align*}

\noindent{\it Case: $\mu>0$. }  By  Lemma~\ref{lm 2.1} with
 $$\tau_1= \frac{2\s_2+N}{2\s_1+ N},\qquad\tau_2=\frac{2\s_2}{2\s_1+ N}$$
 and
 $$r=  \Big(\mu \frac{2\s_2+N}N \Big)^{-\frac{2\s_1+N}{2\s_2}}\frac{(2\s_1+N) M_2}{  M_1 \omega_{_{N-1}}},\qquad d_1=\frac{2\s_2+N}{\omega_{_{N-1}} M_1}\Big(\mu \frac{2\s_2+N}N \Big)^{-\frac{2\s_2+N}{2\s_2}} k>0.$$
Then
 \begin{align*}
  & \Big(\mu \frac{2\s_2+N}N \Big)^{-\frac{2\s_1 +N}{2\s_2}}\frac{(2\s_1 +N) M_2}{  M_1 \omega_{_{N-1}}}
  \\& \geq \Big(\frac{2\s_2+N}{\omega_{_{N-1}} M_1}\Big(\mu \frac{2\s_2+N}N \Big)^{-\frac{2\s_2+N}{2\s_2}} k\Big)^{\frac{2\s_1 +N} {2\s_2+N}}\Big[1-\frac{2s_1+N}{2s_2+N}  \Big(\frac{2\s_2+N}{\omega_{_{N-1}} M_1}\big(\mu \frac{2\s_2+N}N \big)^{-\frac{2\s_2+N}{2\s_2}} k\Big)^{-\frac{2s_2}{2s_2+N}}\Big],
\end{align*}
which is equivalent to
 \begin{align*}
  M_2  \geq  \frac{(2\s_2+N)^{\frac{2\s_1 +N} {2\s_2+N}}}{2\s_1 +N }    (M_1 \omega_{_{N-1}})^{-\frac{2(\s_1 -s_2)} {2\s_2+N}}     k ^{\frac{2\s_1 +N} {2\s_2+N}}\Big[1-\frac{\mu(2s_1+N)}{N(2s_2+N)^\frac{2s_2}{2s_2+N}}(\omega_{N-1}M_1)^\frac{2s_2}{2s_2+N}k^{-\frac{2s_2}{2s_2+N}}\Big],
\end{align*}
that is,
\begin{align*}
 \sum^k_{j=1}\lambda_j (\mu)  & \geq  b_1 |\Omega|^{- \frac{2(s_1-s_2)}{ N}}  k^{1+\frac{2(\s_1-s_2)} {2\s_2+N}}
 \\&\quad -\mu b_1  \frac{2s_1+N}{N(2s_2+N)^{\frac{2s_2}{2s_2+N}}}  \Big((2\pi)^{-N} \omega_{_{N-1}}   c_{s_2} a_0\Big)^{\frac{2s_2}{N+2s_2}}   |\Omega|^{ -\frac{2s_1}{ N}}   k^{1+\frac{2\s_1-4s_2 } {2\s_2+N}}.
\end{align*}
Substituting the value of $a_0$ and $a_{N,s}$ from \eqref{a_0} and \eqref{a_Ns} respectively, we complete the proof.\hfill$
\Box$\medskip

\noindent{\bf Proof of Corollary \ref{cr 1.1}. }  It follows by the nondecreasing monotonicity of $k\mapsto\lambda_k(\mu)$ that
$$\lambda_k(\mu)\geq \frac1k \sum^k_{j=1} \lambda_j(\mu).$$
Applying the above inequality in (\ref{low b}) yields (\ref{low b-k}).\hfill$\Box$

 \setcounter{equation}{0}
\section{Upper bounds}

In order to prove the upper bounds, we need following lemmas.

\begin{lemma}\label{lm 2.1-f}
Let $s\in(0,1)$ and for fixed $z\in \R^N\setminus \{0\}$
\begin{equation}\label{20-11-1}v_z(x)=e^{{\rm i} x\cdot z},\quad\forall\, x\in\R^N,\end{equation}
 then
\begin{equation}\label{2.1}
(-\Delta)^s v_z(x)=   |z|^{2s} v_z(x),\quad\forall\, x\in\R^N.
\end{equation}

\end{lemma}
\noindent{\bf Proof. }  Without loss of generality, it is enough to
prove (\ref{2.1}) with $z=te_1$, where $t>0$ and $e_1=(1,0,\cdots,0)\in\R^N$. For this,
we  write
$$v_t(x) =\mu_z(x_1)=e^{{\rm i} tx_1 },\quad x=(x_1,x')\in \R\times \R^{N-1}.$$
 Using \cite[Lemma 3.1]{CV} we obtian that
  \begin{eqnarray*}
(-\Delta)^s v_t(x)=(-\Delta)^s_{\R} v_t(x_1).
\end{eqnarray*}
 Now we claim  that
\begin{equation}\label{2.2}
(-\Delta)^s_{\R} v_t(x_1)=t^{2s} v_t(x_1),\quad\forall\, x_1\in\R.
\end{equation}
 Indeed, observe that $-\Delta_{\R}v_t:=- (v_t)_{x_1x_1} =t^2 v_t$ in $\R$ and then
$$  (|\xi_1|^{2}-t^2) \widehat{v_t}(\xi_1)=\cF\left(-\Delta_{\R} v_t-t^2v_t\right)(\xi_1)=0, $$
which implies that
$${\rm supp}(\widehat{v_t})\subset \{\pm t\},$$
which in turn implies
$$(|\xi_1|^{2s}-t^{2s}) \widehat{v_t}(\xi_1)=0=\cF\left((-\Delta)^s_{\R} v_t-t^{2s}v_t\right)(\xi_1). $$
and finally
$$\left((-\Delta)^s_{\R} v_t-t^{2s}v_t\right)(\xi_1)=0\quad{\rm in}\ \ \R,$$
which yields
\begin{align*}
(-\Delta)^s v_t(x) &=  (-\Delta)^s_{\R} v_t =  t^{2s}v_t(x),\quad \forall\,x\in\R^N,
 \end{align*}
which completes the proof. \hfill$\Box$\medskip

 Let $\eta_0$ be a $C^2$  increasing function such that ${ \|\eta_0\|_{C^1}},\, \|\eta_0\|_{C^2}\leq 2$,
 $$\eta_0(t)=1\ \ {\rm if} \ \, t\geq 1,\qquad \eta_0(t)=0\ \ {\rm if} \ \, t\leq 0.$$
 For $\sigma>0$, denote
\begin{equation}\label{ws-1}
w_\sigma(x)=\eta_0(\sigma^{-1} \rho(x)),\quad \forall\,  x\in\R^N,
\end{equation}
 where
\begin{equation*}
\rho(x)= \text{dist}(x, \R^N\setminus\Omega) \quad\text{for}\, \, x\in\R^N.
\end{equation*}

Since $\Omega$ is a $C^2$ domain, then $\rho$ is $C^2$ in $\{x\in\R^N: \rho(x)<\delta_0\}$ for
some $\delta_0>0$. Therefore, there is $\sigma_0\in(0,1]$ such that for $\sigma\in(0,\sigma_0]$
$$w_\sigma\in C^2(\R^N).$$
Notice that
$$
w_\sigma \to 1 \quad{\rm in}\  \,\Omega \ \ {\rm as}\ \ \sigma\to0^+.
$$
Moreover, we have that
\begin{equation}\label{test e1}
|\Omega|\geq \int_\Omega w_\sigma dx\geq \int_\Omega w_\sigma^2 dx
\geq    |\Omega_\sigma|,
\end{equation}
thanks to $w_\sigma=1$ in $\Omega_\sigma$.

 \begin{lemma}\label{lm 2.2}
 Let $s\in(0,1)$ and $\Omega$ be a $C^2$ domain,
 then  for $\sigma\in(0,\sigma_0]$
$$|(-\Delta)^s w_\sigma(x)| \leq 2c_{N,s}\omega_{_{N-1}}\sigma^{-2s}  \quad{\rm for}\ \, x\in \Omega.$$
\end{lemma}
\noindent{\bf Proof. }
For $x\in \Omega$, we have that
\begin{align*}
|2w_\sigma(x)-w_\sigma(x+\zeta)-w_\sigma(x-\zeta)|&\leq \min\{  4 , \|w_\sigma\|_{C^2} |\zeta|^2\}
\\&\leq \min\{4, \sigma^{-2}\|\eta_0\|_{C^2} |\zeta|^2\}.
\end{align*}
We use an equivalent definition
\begin{align*}
\frac{2}{c_{N,s}}| (-\Delta)^s w_\sigma(x)| &= \Big|\int_{\R^N} \frac{2w_\sigma(x)-w_\sigma(x+\zeta)-w_\sigma(x-\zeta) }{|\zeta|^{N+2s}} d\zeta\Big|
\\&\leq \int_{\R^N}\frac{\min\{4, \sigma^{-2}\|\eta_0\|_{C^2} |\zeta|^2\}}{|\zeta|^{N+2s}} d\zeta
\\&\leq  2\sigma^{-2} \int_{B_\sigma }\frac{   |\zeta|^2 }{|\zeta|^{N+2s}} d\zeta
+\int_{\R^N\setminus B_\sigma}\frac{ 4}{|\zeta|^{N+2s}} d\zeta
\\&\leq 4\omega_{_{N-1}}\sigma^{-2s},
\end{align*}
where $\|\eta_0\|_{C^2}\leq 2$.
This completes the proof. \hfill$\Box$\medskip

Note that if $w_\sigma$ and $v_z$ are defined by \eqref{20-11-1} and \eqref{ws-1} respectively then
$$(-\Delta)^s(w_\sigma v_z)(z)=v_z(x) (-\Delta)^s w_\sigma(x) +w_\sigma(x) (-\Delta)^s v_z(x)+\cL^s_zw_\sigma(x), $$
where
\begin{equation}\label{20-11-2}\cL^s_z w_\sigma(x)= c_{N,s} \int_{\R^N}\frac{(w_\sigma(x)-w_\sigma(\tilde x))(e^{i\tilde x\cdot z}-e^{i x\cdot z}) }{|x-\tilde x|^{N+2s}}d\tilde x.\end{equation}

\begin{lemma}\label{lm 2.3}
  Let $s\in(0,1)$,  $\Omega$ be a $C^2$ domain
  and $R\geq 1$ be such that $\Omega\subset B_R(0)$, then $\sigma\in(0,\sigma_0]$, $x\in \Omega$ and $|z|>1$ \\[1mm]
$(i)$ for $s\in (\frac12,1)$,
\begin{align*}
\frac{1}{c_{N,s}}| \cL^s_z w_\sigma(x) | \leq\frac{\omega_{_{N-1}}}{1-s}\sigma^{-1}|z|^{2s-1}+\frac{ 4\sigma^{-1}\omega_{_{N-1}}}{2s-1}|z|^{2s-1}+\frac{\omega_{_{N-1}}}{2s}   R ^{-2s};
\end{align*}
$(ii)$ for $s= \frac12$,
\begin{align*}
\frac{1}{c_{N,s}}| \cL^s_z w_\sigma(x) | \leq\frac{\omega_{_{N-1}}}{1-s}\sigma^{-1} |z|^{2s-1}+  4\sigma^{-1} \omega_{_{N-1}} (\log |z| +\log (4R))+\frac{\omega_{_{N-1}}}{2s}   R^{-1};
\end{align*}
$(iii)$  for $s\in (0,\frac12)$,
\begin{align*}
\frac{1}{c_{N,s}}| \cL^s_z w_\sigma(x) | \leq\frac{\omega_{_{N-1}}}{1-s}\sigma^{-1}|z|^{2s-1}+ 4\sigma^{-1}\frac{\omega_{_{N-1}}}{1-2s}(4R)^{1-2s}+\frac{\omega_{_{N-1}}}{2s}  R^{-2s}.
\end{align*}

\end{lemma}
\noindent{\bf Proof. }
Note that
$$|e^{i\tilde x\cdot z}-e^{i x\cdot z}|\leq \min\{2, |z|  |\tilde x-x|\}$$
and
$$|w_\sigma(x)-w_\sigma(\tilde x)|\leq \frac2\sigma |x-\tilde x|,\quad |\tilde x|<3R.$$
 For $x\in \Omega$ and $|z|>1$, we have that
\begin{align*}
\frac{1}{c_{N,s}}| \cL^s_z w_\sigma(x) | &\leq   \int_{\R^N} \frac{|w_\sigma(x)-w_\sigma(\tilde x)|\, |e^{i\tilde x\cdot z}-e^{i x\cdot z}| }{|x-\tilde x|^{N+2s}}d\tilde x
\\& \leq \int_{B_{4R}}  \frac{2\sigma^{-1}|x-\tilde x|\,\min\{2, |z|  |\tilde x-x|\} }{|x-\tilde x|^{N+2s}}d\tilde x +\int_{\R^N\setminus B_{4R}} \frac{2 }{|x-\tilde x|^{N+2s}}d\tilde x
\\&\leq 2\sigma^{-1}|z|\int_{B_{\frac1{|z|}}(x)}   |x-\tilde x|^{2-N-2s}d\tilde x +4\sigma^{-1}\int_{B_{4R}\setminus B_{\frac1{|z|}}(x)}   |x-\tilde x|^{1-N-2s}d\tilde x
\\&\quad +\int_{\R^N\setminus B_{R}} \frac{2 }{|\tilde x|^{N+2s}}d\tilde x,
\end{align*}
where
$$2\sigma^{-1}|z|\int_{B_{\frac1{|z|}}(x)}   |x-\tilde x|^{2-N-2s}d\tilde x\leq \frac{\sigma^{-1}\omega_{_{N-1}}}{1-s} |z|^{2s-1},$$
$$\int_{\R^N\setminus B_{R}} \frac{2 }{|\tilde x|^{N+2s}}d\tilde x\leq \frac{\omega_{_{N-1}}}{2s}   R ^{-2s}$$
and
$$4\sigma^{-1}\int_{B_{4R}\setminus B_{\frac1{|z|}}(x)}   |x-\tilde x|^{1-N-2s}d\tilde x  \leq \left\{ \arraycolsep=1pt
\begin{array}{lll}
4\frac{\sigma^{-1}\omega_{_{N-1}}}{2s-1}|z|^{2s-1} \qquad \  &{\rm if }\ \,   s\in(\frac12,1),\\[2mm]
 \phantom{   }
4\sigma^{-1}\omega_{_{N-1}} (\log |z| +\log (4R))   \quad\ \ &{\rm{if}}\  \,s=\frac12,\\[2mm]
 \phantom{   }
4\frac{\sigma^{-1}\omega_{_{N-1}}}{1-2s}(4R)^{1-2s}  \qquad \ &{\rm{if}}\  \,s\in(0,\frac12).
 \end{array}
\right.$$
This completes the proof. \hfill$\Box$\bigskip

\noindent{\bf Proof of Theorem \ref{teo upp}. } Denote
 $$\Psi_k(x,y)=\sum^k_{j=1} \tilde  \phi_j(x)  \phi_j(y) $$
 and
 $$\cF_x(\Psi_k)(z,y)=(2\pi)^{-\frac N2} \int_{\R^N} \Psi_k(x,y)e^{ix\cdot z}dx.$$

  Let
  $v_{\sigma}(z,y) $ be the solution of
\begin{equation}\label{eq 1.1--}
\left\{ \arraycolsep=1pt
\begin{array}{lll}
 (-\Delta)^{\frac{s_2}2} u= w_\sigma e^{{\rm i}x\cdot z}\quad \  &{\rm in}\quad   \Omega,\\[2mm]
 \phantom{ (-\Delta)^{\frac{s_2}2}  }
  u=0\quad \ &{\rm{in}}\  \quad \R^N\setminus \Omega,
\end{array}
\right.
\end{equation}

{\it We claim that  $v_{\sigma}\in H^s_0(\Omega)$ for any $s\in(0,\frac{1+s_2}{2} )$.}

In fact, from \cite[Corollary 1.6]{RS}, for $\beta \in [\frac{s_2}{2}, \frac{1+s_2}{2})$there exists  $C_\beta>0$ such that
$$[v_{\sigma}]_{C^\beta(\Omega_t)} \leq C_\beta t^{ s_2/2 -\beta},\quad\forall\, t\in (0,r_0)$$
for some $r_0>0$,
where $$\Omega_t=\{x\in\Omega:\, \rho (x)>t\}.$$
 From \cite[Proposition 1.1]{RS}, we also have $v_\sigma\in   C^{s_2/2} (\R^N)$.  First we note that for $\beta \in [\frac{s_2}{2}, \frac{1+s_2}{2})$, it holds
\begin{equation}\label{23-11-1}\big| v_{\sigma}(x)-v_{\sigma}(y)\big|\leq C_\beta\max\big\{\rho^{ s_2/2-\beta}(x), \rho^{ s_2/2-\beta}(y)\big\}\big|x-y\big|^{\beta} \quad{\rm for\ any }\ \ x,\, y\in \Omega.\end{equation}
Indeed, to see the above estimate, note that if $|x-y|\leq 2\rho(x)$ then  $\rho(y)\leq 3\rho(x)$, i.e., $x,y\in\Omega_{\frac{\rho(y)}{3}}$ and thus
$$\big| v_{\sigma}(x)-v_{\sigma}(y)\big|\leq C\rho^{ s_2/2-\beta}(y)\big|x-y\big|^{\beta}.$$
Similarly if $|x-y|\leq 2\rho(y)$ then it follows
$$\big| v_{\sigma}(x)-v_{\sigma}(y)\big|\leq C\rho^{ s_2/2-\beta}(x)\big|x-y\big|^{\beta}.$$
Finally when $|x-y|>2\max\{\rho(x), \rho(y)\}$, then from \cite{RS} we have
$$\big| v_{\sigma}(x)-v_{\sigma}(y)\big|\leq |v_\sigma(x)|+|v_\sigma(y)|\leq C(\rho^\frac{s_2}{2}(x)+\rho^\frac{s_2}{2}(y))\leq C\max\big\{\rho^{ s_2/2-\beta}(x), \rho^{ s_2/2-\beta}(y)\big\}\big|x-y\big|^{\beta}.$$
Hence \eqref{23-11-1} follows.

 Next, for any given $s\in(\frac{s_2}{2}, \frac{1+s_2}{2})$. We fix some $\beta>0$ such that $\beta\in(s, \frac{1+s_2}{2})$. Then using \eqref{23-11-1}, we have that

\begin{align*}\int_{\Omega} \int_{ \Omega} \frac{|v_\sigma(x)-v_\sigma(y)|^2}{|x-y|^{N+2s}}  dy dx &\leq C\int_{\Omega} \int_{ \Omega}  \bigg(\rho^{s_2-2\beta}(x)+\rho^{s_2-2\beta}(y)\bigg)|x-y|^{2\beta-N-2s}dydx\\
&\leq 2\int_{\Omega} \int_{ \Omega} \rho^{s_2-2\beta}(y)|x-y|^{2\beta-N-2s}dy dx\\
&\leq   2\int_{\Omega}\bigg( \int_{|x-y|\leq\rho(y)}|x-y|^{2\beta-N-2s}dx\bigg)\rho^{s_2-2\beta}(y) dy\\
&\quad\ \ +2\int_{\Omega}\bigg(\int_{|x-y|>\rho(y)}|x-y|^{2\beta-N-2s}dx\bigg)\rho^{s_2-2\beta}(y) dy\\
&\leq C\int_{\Omega} \rho^{s_2-2\beta}(y) dy+ C\int_{\Omega}\rho^{2\beta-2s+s_2-2\beta}(y)dy\\
&\leq C\int_{t=0}^{r_0}\int_{\rho(y)=t}t^{s_2-2\beta}dSdt
+C\int_{t=0}^{r_0}\int_{\rho(y)=t}t^{s_2-2s}dSdt\\
&<\infty,
\end{align*}
for  our choice that $\frac{s_2}{2}<s <\beta<\frac{1+s_2}{2}$.
Similarly, we also show that
\begin{align*}\int_{\Omega^c} \int_{ \Omega} \frac{|v_\sigma(x)-v_\sigma(y)|^2}{|x-y|^{N+2s}}  dy dx &= \int_{\Omega^c} \int_{ \Omega} \frac{|v_\sigma(y)|^2}{|x-y|^{N+2s}}  dy dx\leq \int_{\Omega^c} \int_{ \Omega} \rho(y)^{2s_2}|x-y|^{-(N+2s)}dy dx<\infty,
\end{align*}
as $s<\frac{1+s_2}{2}$ and for any $y\in\Omega$, it follows $\{x:\in\Omega^c: |x-y|<\rho(y)\}=\emptyset$.
 Hence, $v_\sigma \in H^s_0(\Omega)$ for $s<\frac{1+s_2}2 $.

Denote
$$\tilde v_{\sigma}(x,z)=  (-\Delta)^{\frac{s_2}2}   v_{\sigma}(x,z) $$
and
 $$v_{\sigma,k}(z,y)  := v_{\sigma}(y,z)- \sum^k_{j=1} \big\langle \tilde v_{\sigma}(.,z), \tilde\phi_j\big\rangle_{L^2( \R^N)} \phi_j(y).  $$

Note that
$$(-\Delta)^{\frac{s_2}2} _y  v_{\sigma,k}(z,y)  ={ \tilde v_{\sigma}(y,z)}-\sum^k_{j=1}   \big\langle \tilde v_{\sigma}(.,z), \tilde\phi_j\big\rangle_{L^2({ \R^N})} \tilde \phi_j(y), $$
thus we  get $v_{\sigma,k}(z,\cdot) \in \cH_{0,k+1}(\Omega)$
and
the Rayleigh-Ritz formula shows that
\begin{align*}
\lambda_{k+1}  \Big| \int_{ \R^N}(-\Delta)^{\frac{s_2}2}_y  v_{\sigma,k}(z,y)   dy \Big|^2 \leq \Big| \int_{ \R^N} (-\Delta)^\frac{s_1}{2}_{y}  v_{\sigma,k}(z,y) dy \Big|^2
\end{align*}
for any $z\in\R^N$ and $\sigma>0$.
Thus we can conclude that
$$\lambda_{k+1}(\mu)\leq \inf_{\sigma>0} \frac{\int_{B_r} \int_{ \R^N}   (-\Delta)^{s_1}_{y}  v_{\sigma,k}(z,y) \,\overline{ v_{\sigma,k}(z,y)} dydz }{\int_{B_r} \int_{ \R^N}  (-\Delta)^{s_2}_y  v_{\sigma,k}(z,y)  \, \overline{v_{\sigma,k}(z,y)} dy dz }.$$

An elementary calulation yields that
\begin{align*}
&\int_{B_r} \int_{ \R^N}  (-\Delta)^{s_2}_y  v_{\sigma,k}(z,y)  \, \overline{v_{\sigma,k}(z,y)} dy dz
\\ &= \int_{B_r} \int_{\R^N}\big| (-\Delta)^{\frac{s_2}2}_y  v_{\sigma,k}(z,y)   \big|^2dydz
\\ &= \int_{B_r} \int_{\R^N}\big| (-\Delta)^{\frac{s_2}2}_y  v_{\sigma}(z,y)   \big|^2dy dz-\int_{B_r} \sum^k_{j=1} \Big|\big\langle \tilde v_{\sigma}(.,z), \tilde\phi_j\big\rangle_{L^2( \R^N)}\Big|^2\int_\Omega \tilde \phi_j^2  dxdz
  \\&=  \int_{B_r} \int_{{\R^N}}  w_\sigma ^2      dxdz
 -   \sum^k_{j=1} \int_{B_r}  \Big|\big\langle \tilde v_{\sigma}(.,z), \tilde\phi_j\big\rangle_{L^2( \R^N)}\Big|^2dz
 \end{align*}
 and
\begin{align*}
& \int_{B_r}\Big( \int_{{\R^N}}  (-\Delta)^{s_1}_{y}v_{\sigma,k}(z,y)  \, \overline{ v_{\sigma,k}(z,y)}dy dz
\\ &= \int_{B_r} \int_{\R^N} \Big( (-\Delta)^{s_1}_{y}  v_{\sigma}( y,z)\,  \overline{ v_{\sigma}( y,z) }  \Big)  dydz-\int_{B_r} \sum^k_{j=1} \Big|\big\langle \tilde v_{\sigma}(.,z), \tilde\phi_j\big\rangle_{L^2( \R^N)}\Big|^2\int_\Omega  \phi_j (-\Delta)^{s_1}\phi_j   dxdz
 \\&= \int_{B_r} \int_{\R^N}\Big( (-\Delta)^{s_1-s_2}_{y} (w_\sigma e^{{\rm i}x\cdot z}) \Big) \,\overline{(w_\sigma e^{{\rm i}x\cdot z}) }   dydz
   - \sum^k_{j=1}\lambda_j(\mu)\int_{B_r}   \Big|\big\langle \tilde v_{\sigma}(.,z), \tilde\phi_j\big\rangle_{L^2( \R^N)}\Big|^2dz
  \\&=  \int_{B_r} \int_{\R^N} \Big(  \cL_z^{s_1-s_2}(w_\sigma)\, \overline{w_\sigma e^{{\rm i}x\cdot z}}+w_\sigma  (-\Delta)^{s_1-s_2}w_\sigma    +  |z|^{2(s_1-s_2)}   w_\sigma^2(x) \Big)       dxdz
\\&\qquad -   \sum^k_{j=1} \lambda_j(\mu)\int_{B_r}  \Big|\big\langle \tilde v_{\sigma}(.,z), \tilde\phi_j\big\rangle_{L^2( \R^N)}\Big|^2dz.
 \end{align*}

Since $\Omega$ is $C^2$, there exists $t_0\in(0,1)$ such that
$$|\partial \Omega_t|\leq 2 |\partial\Omega|\quad{\rm for }\ t\in(0,t_0].$$
From (\ref{test e1}),
\begin{align*} |\Omega|\geq \int_{\Omega} w_\sigma^2(y)dy>  |\Omega_\sigma|
\geq |\Omega|-2\sigma |\partial\Omega|>\frac{|\Omega|}2,
\end{align*}
when $\sigma$ is small enough. We choose $r>r_0$ for some $r_0>1$ such that $r^{-\frac {s_1-s_2}{2}}=\sigma$ satisfies the above relation and $\sigma< \frac1{4|\partial\Omega|}|\Omega|$. Therefore,
 \begin{equation}\label{sigma-1}
 |\Omega|\geq \int_{\Omega} w_\sigma^2(y)dy \geq |\Omega|-2 r^{- \frac {s_1-s_2}{2}} |\partial\Omega|>\frac{|\Omega|}2.
  \end{equation}

Note that for  $s=s_1-s_2$ and $r>r_0$
 \begin{align*}
  &\frac1{c_{N,s}} \int_{B_r} \int_{\R^N} |\cL^{s}_z w_\sigma \, \overline{w_\sigma e^{{\rm i}x\cdot z}}|    dydz
\\ &\leq  \frac{\omega_{_{N-1}}r^{N+2s-1}}{\sigma(N+2s-1)(1-s)}|\Omega|
 +\frac1\sigma \varphi_{s}(r,R) |\Omega| +  \frac{\omega_{_{N-1}}^2r^{N}}{2sN}R^{-2s}|\Omega|
 \\& \leq   \frac{\omega_{_{N-1}}|\Omega| }{ (N+2s-1)(1-s)}r^{N+ 2s+ \frac {s}{2}-1 } +r^{ \frac {s}{2}} \varphi_s(r,R) |\Omega| +\frac{\omega_{_{N-1}}^2r^{N}}{  2sN}R^{-2s}|\Omega|
 \\& \leq c_1(N,s_1, s_2)|\Omega|\, r^{N+\max\{ 2s+\frac {s}{2}-1,\,  \frac {s}{2}\}},
 \end{align*}
 where

\begin{equation}\label{3.1"}
\varphi_{s}(r,R) =\left\{ \arraycolsep=1pt
\begin{array}{lll}
\frac{ 4\omega_{_{N-1}}r^{N+2s-1}}{(2s-1)(N+2s-1)} \qquad \  &{\rm if }\quad   s\in(\frac12,1),\\[2mm]
 \phantom{   }
 \frac{ 4\omega_{_{N-1}}r^{N}(\log r+\log (4R)) }{N}  \qquad \ &{\rm{if}}\  \,s=\frac12,\\[2mm]
 \phantom{   }
\frac{ 4\omega_{_{N-1}}^2r^{N}  }{N(1-2s)}(4R)^{1-2s}  \qquad \ &{\rm{if}}\  \,s\in(0,\frac12).
 \end{array}
\right.
\end{equation}
Moreover, we have that
\begin{align*}
  \int_{B_r}  \int_{\R^N} |w_\sigma  (-\Delta)^{s}w_\sigma| dy dz&\leq   2c_{N,s} \frac{\omega_{_{N-1}}^2r^{N}  }{N}  \sigma^{-2s} |\Omega|= c_2(N,s_1, s_2)|\Omega|r^{N+  s^2}.
 \end{align*}
Furthermore,
\begin{align*}
 \int_{B_r}  \int_{\R^N}  |z|^{2s}   w_\sigma^2(x)  dy dz=  \frac{\omega_{_{N-1}}r^{N+2s}  } {N+2s}\int_{\Omega} w_\sigma^2(y) dy.
 \end{align*}
Let

 $$\delta_0=\max\big\{ 2s+ \frac {s}{2}-1,\,  \frac {s}{2}, s^2 \big\}=  \frac {s}{2},$$
 where we have used the hypothesis that $s_1<\frac{1+s_2}{2}$.

%
%
Let
  \begin{align*}
   P_j:&=\int_{B_r}\Big|\big\langle \tilde v_{\sigma}(.,z), \tilde\phi_j\big\rangle_{L^2( \R^N)}\Big|^2dz
   \\&\leq\int_{\R^N}\Big| \int_{\R^N} w_\sigma e^{{\rm i}x\cdot z} \tilde\phi_j\Big|^2dz
 \ \leq(2\pi)^N \!\int_{\R^N}\! \Big| \cF(w_\sigma \tilde\phi_j)\Big|^2dz
 \\&\qquad\qquad \qquad \qquad \qquad\qquad    =\ \, (2\pi)^N  \int_{\R^N}  (w_\sigma \tilde\phi_j)^2 dx\ \,   \leq\, (2\pi)^N,
  \end{align*}
using the Parseval's inequality.   Note that
  \begin{align*}
 \lambda_{k+1}(0) \leq  \frac{\frac{ \omega_{_{N-1}}}{N+2(s_1-s_2)} r^{N+2(s_1-s_2)}\displaystyle\int_{\Omega} w_\sigma^2(y)dy +O(1)|\Omega| r^{N+\frac{s}{2}}  -\displaystyle  \sum^k_{j=1}\lambda_j(0) P_j   }{\frac{\omega_{_{N-1}}}{N } r^{N } \displaystyle \int_{\Omega} w_\sigma^2(y) dy -   \displaystyle  \sum^k_{j=1}  P_j   },
 \end{align*}
 and now taking
 $$Q_1=\frac{ \omega_{_{N-1}}}{N+ 2(s_1-s_2)} r^{N+2(s_1-s_2)}\int_{\Omega} w_\sigma^2(y)dy +O(1)|\Omega| r^{N+\frac{s}{2}}\quad{\rm and}\quad  Q_2= \frac{\omega_{_{N-1}}}{N } r^{N }   \int_{\Omega} w_\sigma^2(y) dy, $$
we have that
 \begin{align*}
0&\leq \frac{Q_1- \sum^k_{j=1}\lambda_j(0) P_j }{Q_2 - \sum^k_{j=1}  P_j }-\lambda_{k+1}(0)
\\[2mm]&= \frac{\big(Q_1-Q_2\lambda_{k+1}(0)\big)+ \sum^k_{j=1}\Big(\lambda_{k+1}(0)-\lambda_j(0)\Big) P_j }{Q_2 - \sum^k_{j=1}  P_j  }
\\[2mm]&\leq \frac{\big(Q_1-Q_2\lambda_{k+1}(0)\big)+(2\pi)^N\sum^k_{j=1}\Big(\lambda_{k+1}(0)-\lambda_j(0)\Big) }{Q_2 -(2\pi)^N k },
\end{align*}
since $\lambda_{k+1}(0)\geq \lambda_{j}(\mu)$ for $j< k+1$ and $ P_j \in \big(0, (2\pi)^N]$.
Therefore
$$
 0<\lambda_{k+1}(0) \leq \frac{\frac{ \omega_{_{N-1}}}{N+2(s_1-s_2)} r^{N+2(s_1-s_2)}\displaystyle\int_{\Omega} w_\sigma^2(y)dy +O(1)|\Omega| r^{N+\frac{s_1-s_2}{2}}  -(2\pi)^N\displaystyle  \sum^k_{j=1}\lambda_j(0)}
   {\frac{\omega_{_{N-1}}}{N } r^{N }\displaystyle\int_{\Omega} w_\sigma^2(y) dy  -(2\pi)^N k}.
$$

There exists $k_0\geq 1$ such that for $k\geq k_0$, we can choose $r>r_0 $ satisfying
\begin{equation}\label{23-11-3}\frac{\omega_{_{N-1}}}{N } r^{N } \int_{\Omega} w_\sigma^2(y) dy=(2\pi)^N (k+1).
\end{equation}

As a consequence, using  (\ref{sigma-1}) and the above relation, we obtain   for $k\geq k_0$
\begin{align*}
\sum^k_{j=1}\lambda_j(0)& \leq(2\pi)^{-N} \frac{ \omega_{_{N-1}}}{N+2(s_1-s_2)} r^{N+2(s_1-s_2)}\int_{\Omega} w_\sigma^2(y)dy +c |\Omega|  r^{N+\frac{s_1-s_2}{2}}
\\&= (2\pi)^{2(s_1-s_2)} \frac{N}{N+2(s_1-s_2)} \Big(N^{-1} \omega_{_{N-1}}\int_{\Omega} w_\sigma^2(y)dy\Big)^{-\frac{2(s_1-s_2)}{N}} (k+1)^{1+2\frac{\s_1-s_2} { N}}
\\& \qquad +c_3 k^{1+\frac12\frac{s_1-s_2}{N}}
\\&\leq b_3 |\Omega|^{-\frac{2(s_1-s_2)}{ N}}  k^{1+2\frac{\s_1-s_2} { N}} +  c_4  k^{ 1+\frac32\frac{s_1-s_2}{N}  }+c_3  k^{1+\frac12\frac{s_1-s_2}{N}} ,
\end{align*}
 where $c_3, c_4>0$ depends on $N, s_1,s_2$, $\Omega$,
$$ (k+1)^{1+2\frac{\s_1-s_2} { N}}\leq k ^{1+2\frac{\s_1-s_2} { N}}+c_5k^{2\frac{\s_1-s_2} { N}}\quad{\rm for}\ \, k\geq 1$$ and
  we use the fact that for $c_7,c_8>0$
\begin{align*}
\Big(\int_{\Omega} w_\sigma^2(y)dy\Big)^{-\frac{2(s_1-s_2)}{N} } &\leq \Big(|\Omega|-2 r^{- \frac {s_1-s_2}{2}} |\partial\Omega|\Big)^{-\frac{2(s_1-s_2)}{N} }
\\& \leq   |\Omega|^{-\frac{2(s_1-s_2)}{N} }+c_6 r^{- \frac {s_1-s_2}{2}}
\\& \leq   |\Omega|^{-\frac{2(s_1-s_2)}{N} } +c_7 k^{- \frac {s_1-s_2}{2N}}.
\end{align*}
Here for $k\leq k_0$, we only have to adjust the constant $c_3$ or $c_4$.     \hfill$\Box$

\bigskip\medskip

\noindent{\small {\bf Acknowledgements:}  This work of H.~Chen is supported by NNSF of China, No: 12071189 and 12001252 and the Jiangxi Provincial Natural Science Foundation, No: 20202BAB201005, 20202ACBL201001.

The research of M.~Bhakta is partially supported by the  SERB MATRICS grant MTR/2017/000168 and SERB WEA grant WEA/2020/000005.


\begin{thebibliography}{99}

 
\bibitem {BCP} M. Bhakta, S. Chakraborty and P. Pucci,  Fractional Hardy-Sobolev equations with nonhomogeneous terms, {\it preprint }, arXiv:2008.01118.

\bibitem {BM} M. Bhakta and D. Mukherjee, Multiplicity results for $(p,q)$ fractional elliptic equations involving critical nonlinearities, {\it Adv. Diff. Equat. 24(3-4)}, 185-228 (2019). 

\bibitem{BM-2}  M. Bhakta and D. Mukherjee, Multiplicity results and sign changing solutions of non-local equations with concave-convex nonlinearities. {\it Differential Integral Equations} 30 (2017), no. 5-6, 387--422. 


\bibitem {B-1}  S. Biagi, S. Dipierro, E. Valdinoci and E. Vecchi,
Mixed local and nonlocal elliptic operators: regularity and maximum principles,
{\it arxiv.org/pdf/2005.06907. }

\bibitem {CS0} L. Caffarelli, S. Salsa and L. Silvestre,
Regularity estimates for the solution and the free boundary to the
obstacle problem for the fractional Laplacian, {\it Invent.
Math. 171,} 425-461 (2008).


\bibitem {C} H-Y. Chen, Liouville theorem for the fractional Lane-Emden equation in an unbounded domain, \emph{J. Math. Pures et Appl. 111}, 21-46 (2018).


\bibitem {CFQ} H-Y. Chen, P. Felmer and A. Quaas, Large solution to elliptic  equations involving fractional Laplacian, {\it  Ann.  Inst. Henri Poincar\'{e}-AN 32,} 1199-1228 (2015).

  \bibitem {CQL} H. Chen,  R. Qiao,  P.  Luo and  D. Xiao, Lower and upper bounds of Dirichlet eigenvalues for totally characteristic degenerate elliptic operators, {\it Sci. China Math. 57(11),} 2235-2246 (2014).

\bibitem {CQ} H-Y. Chen and A. Quaas, Classification of isolated singularities of nonnegative solutions to fractional semilinear elliptic equations and the existence results,
\emph{J. London Math. Soc. 97(2),} 196-221 (2018).


\bibitem {CV} H-Y. Chen and L. V\'eron, Initial trace of positive solutions to fractional diffusion equations with absorption, \emph{J. Funct. Anal. 276,} 1145-1200 (2019).

\bibitem {CT} H-Y. Chen and T. Weth, The Dirichlet Problem for the Logarithmic Laplacian, \emph{Comm. Part. Diff. Eq. 44,} 1100-1139 (2019).

\bibitem {CZ}  H. Chen and A. Zeng, Universal inequality and upper bounds of eigenvalues for
   non-integer poly-Laplacian on a bounded domain, {\it  Calc.Var. Part. Diff. Eq.  56(5) } (2017).


     \bibitem {CgY} Q.  Cheng and H. Yang,  Bounds on eigenvalues of Dirichlet Laplacian, {\it   Math. Ann. 337,} 159-175   (2007).


    \bibitem {CgW} Q. Cheng and G. Wei,  A lower bound for eigenvalues of a clamped plate problem, {\it  Calc.Var. Part. Diff. Eq. 42(3/4),} 579-590 (2011).

 \bibitem {EGE} E. Di Nezza, G. Palatucci and E. Valdinoci,  Hitchhiker's guide to the fractional Sobolev spaces,
{\it  Bull.   Sci. Math.   136(5)},  521-573 (2012).


\bibitem {E-3} E. Elshahed, A fractional calculus model in semilunar heart valve vibrations,{\it International Mathematica symposium} (2003).

\bibitem{FKV13}
M.~Felsinger, M.~Kassmann and P.~Voigt, The Dirichlet problem for nonlocal operators,  {\it Math. Zeit. 279,} 779--809 (2015).

\bibitem{F1}  R. Frank,  Eigenvalue bounds for Schr\"odinger operators with complex potentials. III. {\it Trans. Amer. Math. Soc. 370(1)}, 219-240 (2018).

\bibitem{F2}  R. Frank,  Eigenvalue bounds for the fractional Laplacian: a review. {\it Recent developments in nonlocal theory}, 210--235, De Gruyter, Berlin, (2018).

\bibitem{Frank} R.  Frank,  E. Lenzmann and L. Silvestre, Uniqueness of radial solutions for the fractional Laplacian, {\it Comm. Pure Appl. Math. 69(9)}, 1671-1726  (2016).

\bibitem {G} L. Geisinger,  A short proof of Weyl's law for fractional differential operators,   {\it J. Math. Phys. } 011504, Doi 10.1063/1.4861935 (2014).

\bibitem{GS1} S.  Goyal and K. Sreenadh,  On the Fu\v{c}ik spectrum of non-local elliptic operators,
{\it  Nonlinear Diff. Eq.    Appl.   21(4),}  567-588 (2014).


 \bibitem {HY} E. M. Harrell II and S. Y. Yolcu,  Eigenvalue inequalities for Klein-Gordon operators,  {\it  J. Funct. Anal. 256(12),} 3977-3995 (2009).

\bibitem{Hajaiej1} H. Hajaiej, Existence of minimizers of functional involving the fractional gradient in the absence of compactness, symmetry and monotonicity, {\it Journal of Mathematical Analysis and Appl. 399(1),} 17-26 (2013).

\bibitem{Hajaiej2} H. Hajaiej, On the optimality of the conditions used to prove the symmetry of the minimizers of some fractional constrained variational problems, {\it Annales de l’Institut Henri Poincare. 14(5),} 1425-1433 (2013).

\bibitem{Hajaiej3} H. Hajaiej, Symmetry of minimizers of some fractional problems, {\it Applicable Analysis 94(4),}: 1-7 (2014).

\bibitem {JXL}  T. Jin, Y. Li and J. Xiong,  On a fractional Nirenberg problem, part I: Blow up analysis and compactness of solutions, {\it J. Eur. Math. Soc. 16(6),}   1111-1171 (2014).




\bibitem {KM} M. Kassmann and A. Mimica, Intrinsic scaling properties for nonlocal operators, {\it J. Eur. Math. Soc. 19(4)}, 983-1011 (2013).

\bibitem {K} P. Kr\"oger,  Estimates for sums of eigenvalues of the Laplacian, {\it J. Funct. Anal. 126(1),} 217-227 (1994).

\bibitem {L}  A. Laptev, Dirichlet and Neumann eigenvalue problems on domains in Euclidean spaces, {\it J. Funct. Anal. 151(2)} 531-545 (1997).

\bibitem {LY} P. Li and  S.-T.Yau,   On the Schr\"odinger equation and the eigenvalue problem. {\it Commun. Math. Phys. 88(3),} 309-318 (1983).

  \bibitem {L} E. Lieb,   The number of bound states of one-body Schr\"odinger operators and the Weyl problem, {\it Proc. Sym. Pure Math. 36,} 241-252 (1980).

  \bibitem {M-4}  R. L. Magin, Fractional calculus in bioenginering 1, 2, 3,
{\it Critical Reviews in Biomedical Engineering, } 1-1377, 32 (2004).

\bibitem {M-5}  R.L. Magin, S. Boregowda, and C. Deodhar, Modelling of pulsating peripheral bioheat transfusing fractional calculus and constructal theory,
{\it Journal of Design $\&$ Nature 1,}   18-33 (2007).

\bibitem {M-6}   R.L. Magin and M. Ovadia, Modeling the cardiac tissue electrode interface using fractional calculus, {\it Journal of Vibration and Control 19},1431-1442 (2009)

 \bibitem {M}  A. Melas,  A lower bound for sums of eigenvalues of the Laplacian, {\it Proc. Am. Math. Soc. 131(2),} 631-636 (2003).

\bibitem{musina-nazarov} R. Musina and A.I. Nazarov, On fractional Laplacians. {\it Comm. Part. Diff. Eq.  39,} 1780-1790  (2014).

 \bibitem {P} G. P\'olya,   On the Eigenvalues of Vibrating Membranes
 (In Memoriam Hermann Weyl), {\it  Proc. Lond. Math. Soc. 3(1),} 419-433 (1961).

\bibitem {RS} X. Ros-Oton and J. Serra, The Dirichlet problem for the fractional
laplacian: regularity up to the boundary, {\it J. Math. Pures Appl. 101},  275-302 (2014).

\bibitem {RS1} X. Ros-Oton and J. Serra, The Pohozaev identity for the fractional Laplacian,
{\it Arch.   Ration. Mech.   Anal. 213,} 587-628 (2014).


\bibitem {SV} R. Servadei and E. Valdinoci,    Variational methods for non-local operators of elliptic type, {\it  Discrete Contin. Dyn. Syst. 33(5),} 2105-2137 (2013).

 \bibitem {Tr}  H. Triebel,  Theory of function space, {\it Monographs in Mathematics 78, },  Birkhuser Verlag, Basel  1983.

\bibitem {W} H. Weyl,  Das asymptotische Verteilungsgesetz der Eigenwerte linearer partieller Differentialgleichungen (mit einer Anwendung auf die Theorie der Hohlraumstrahlung). {\it Math. Ann. 71(4),} 441-479 (1912).

\bibitem {YH} Y. Wang and H. Hichem,  Kr\"oger's type upper bounds for Dirichlet eigenvalues
of the fractional laplacian, {\it preprint }   (2020).


\bibitem {YY} S.Y. Yolcu and T. Yolcu,   Estimates for the sums of eigenvalues of the fractional Laplacian on a bounded domain, {\it Commun. Contemp. Math. 15(3),} 1250048 (2013).



\bibitem {YY1}  S.Y. Yolcu and T. Yolcu,  Sharper estimates on the eigenvalues of Dirichlet fractional Laplacian, {\it Discr. Cont. Dyn. Syst. 35(5)},     2209-2225 (2014).









\end{thebibliography}
\end{document}